\documentclass[11pt,leqno]{article}
\usepackage{amsmath,amssymb,amscd,latexsym}
\newcommand{\Co}{\mathrm{Co}}
\newcommand{\Sp}{\mathrm{Sp}}
\newcommand{\Sa}{{\mathbb{S}}}
\newcommand{\SO}{\mathrm{SO}}
\newcommand{\PSO}{\mathrm{PSO}}
\newcommand{\proj}{\mathrm{proj}}
\newcommand{\et}{\mathrm{\acute{e}t}}
\newcommand{\Frob}{\mathrm{Frob}}
\newcommand{\ann}{\mathrm{an}}
\newcommand{\oUh}{\overline{\Uh}}
\newcommand{\oeX}{\overline{\eX}}
\newcommand{\oF}{\overline{F}}
\newcommand{\OF}{\overline{\F}}
\newcommand{\oX}{\overline{X}}
\newcommand{\hzeta}{\hat{\zeta}}
\newcommand{\ox}{\overline{x}}
\newcommand{\C}{{\mathbb{C}}}

\newcommand{\F}{{\mathbb{F}}}

\newcommand{\Q}{{\mathbb{Q}}}
\newcommand{\R}{{\mathbb{R}}}
\newcommand{\Z}{{\mathbb{Z}}}

\newcommand{\ar}{\mathrm{ar}}
\newcommand{\ddet}{\mathrm{det}}
\newcommand{\diss}{\mathrm{dis}}
\newcommand{\id}{\mathrm{id}}

\newcommand{\dyn}{\mathrm{dyn}}
\renewcommand{\mod}{\mathrm{mod}\;}

\newcommand{\ord}{\mathrm{ord}}

\newcommand{\sgn}{\mathrm{sgn}\,}
\newcommand{\spec}{\mathrm{spec}\,}
\newcommand{\Arg}{\mathrm{Arg}\,}

\newcommand{\Fr}{\mathrm{Fr}}

\newcommand{\GL}{\mathrm{GL}\,}
\newcommand{\Gr}{\mathrm{Gr}\,}

\newcommand{\Imm}{\mathrm{Im}\,}
\newcommand{\imm}{\mathrm{im}\,}
\newcommand{\ind}{\mathrm{ind}}
\newcommand{\Ker}{\mathrm{Ker}\,}
\newcommand{\prim}{\mathrm{prim}}
\newcommand{\RRe}{\mathrm{Re}\,}
\newcommand{\ssp}{\mathrm{sp}\,}
\newcommand{\Tr}{\mathrm{Tr}}
\newcommand{\tr}{\mathrm{tr}}
\newcommand{\vol}{\mathrm{vol}}
\newcommand{\Ah}{{\cal A}}
\newcommand{\Ch}{{\cal C}}
\newcommand{\Dh}{{\cal D}}

\newcommand{\Fh}{{\cal F}}

\newcommand{\Hh}{{\mathcal H}}
\newcommand{\Lh}{{\mathcal L}}
\newcommand{\Mh}{{\cal M}}
\newcommand{\Oh}{\mathcal{O}}
\newcommand{\Rh}{{\cal R}}

\newcommand{\Uh}{\mathcal{U}}

\newcommand{\eX}{{\cal X}}
\newcommand{\ea}{\mathfrak{a}}
\newcommand{\eo}{\mathfrak{o}}
\newcommand{\ep}{\mathfrak{p}}
\newcommand{\oH}{\bar{H}}
\newcommand{\tzeta}{\tilde{\zeta}}
\newcommand{\ohne}{\setminus}
\newcommand{\silo}{\stackrel{\sim}{\longrightarrow}}
\newcommand{\tei}{\, | \,}
\newcommand{\ent}{\;\widehat{=}\;}
\newcommand{\hullet}{{\scriptscriptstyle \bullet}}
\newcommand{\verk}{\mbox{\scriptsize $\,\circ\,$}}
\newcommand{\halb}{\frac{1}{2}}

\newcommand{\dis}{\displaystyle}
\newcommand{\Rpp}[2]{\prod_{#1}^{#2} \hspace{-4mm} {\rule[-4pt]{2.5mm}{0.2mm}} \hspace{3mm}}
\newcommand{\RRp}[2]{\prod_{#1}^{#2} \hspace{-7mm} {\rule[-4pt]{2.5mm}{0.2mm}} \hspace{5mm}}
\newcommand{\RRP}[2]{\prod_{#1}^{#2}\hspace{-9.5mm} {\rule[-2pt]{2.5mm}{0.2mm}}\hspace*{7mm}}
\newtheorem{theorem}{Theorem}[section]
\newtheorem{prop}[theorem]{Proposition}

\newtheorem{cor}[theorem]{Corollary}

\newtheorem{remarks}[theorem]{Remarks}
\newtheorem{conj}[theorem]{Conjecture}
\newtheorem{fact}[theorem]{Fact}
\newtheorem{punkt}[theorem]{$\!\!$}
\newenvironment{example}{\bigskip \noindent {\bf Example}}{}
\newenvironment{proof}{\bigskip \noindent {\bf Proof}}{\mbox{}\hspace*{\fill}$\Box$}
\newenvironment{remarknn}{\bigskip \noindent {\bf Remark}}{}
\newenvironment{remarksnn}{\bigskip \noindent {\bf Remarks}}{}
\newenvironment{proofof}{\bigskip \noindent {\bf Proof of}}{\mbox{}\hfill$\Box$}
\begin{document}

\title{Arithmetic Geometry and Analysis on Foliated Spaces}
\author{Christopher Deninger}
\date{\today}
\maketitle

\section{Introduction}

For the arithmetic study of varieties over finite fields powerful cohomological methods are available which in particular shed much light on the nature of the corresponding zeta functions. These investigations culminated in Deligne's proof of an analogue of the Riemann conjecture for such zeta functions. This had been the hardest part of the Weil conjectures. For algebraic schemes over $\spec \Z$ and in particular for the Riemann zeta function no cohomology theory has yet been developed that could serve similar purposes. For a long time it had even been a mystery how such a theory could look like even formally. In these lectures we first describe the shape that a cohomological formalism for algebraic schemes over the integers should take. We then discuss how it would relate to conjectures on arithmetic zeta-functions and indicate a couple of consequences of the formalism that can be proved using standard methods. As it turns out there is a large class of dynamical systems on foliated manifolds whose reduced leafwise cohomology has several of the expected structural properties of the desired cohomology for algebraic schemes. Comparing the arithmetic and dynamical pictures leads to some insight into the basic geometric structures that dynamical systems relevant for $L$-functions of varieties over number fields should have. We also discuss relations between these ideas on dynamical systems and ``Arithmetic Topology''. In that subject one studies analogies between number theory and the theory of $3$-manifolds. 

The present report was written for a lecture series at the Southwestern Center for Arithmetic Algebraic Geometry at the University of Arizona. It is mostly based on the papers \cite{D10}, \cite{D14} and \cite{D15}, properly updated.

I would like to thank Minhyong Kim very much for the invitation to give these lectures.

\section{Geometric zeta- and $L$-functions}

Consider the Dedekind zeta function of a number field $K / \Q$
\[
\zeta_K (s) = \prod_{\ep} (1 - N\ep^{-s})^{-1} = \sum_{\ea} N \ea^{-s} \quad \mbox{for} \; \RRe s > 1 \; .
\]
Here, $\ep$ runs over the prime ideals of the ring of integers $\eo = \eo_K$ of $K$. 
The function $\zeta_K (s)$ has a holomorphic continuation to $\C \setminus \{ 1 \}$ with a simple pole at $s =1$. To its finite Euler factors
\[
\zeta_{\ep} (s) = (1 - N\ep^{-s})^{-1}
\]
we add Euler factors corresponding to the archimedian places $\ep \tei \infty$ of $K$
\[
\zeta_{\ep} (s) = \left\{ \begin{array}{ccl}
 2^{-1/2} \, \pi^{-s/2} \, \Gamma (s / 2) & \mbox{if} & \ep \tei \infty \quad \mbox{is real} \\
(2 \pi)^{-s} \Gamma (s) & \mbox{if} & \ep \tei \infty \quad \mbox{is complex}
\end{array} \right.
\]
The completed zeta function
\[
\hat{\zeta}_K (s) = \zeta_K (s) \prod_{\ep \tei \infty} \zeta_{\ep} (s)
\]
is holomorphic in $\C \setminus \{ 0,1 \}$ with simple poles at $s = 0,1$ and satisfies the functional equation:
\[
\hat{\zeta}_K (1-s) = |d_K|^{s- 1/2} \hat{\zeta}_K (s) \; .
\]
Here $d_K$ is the discriminant of $K$ over $\Q$. The zeroes of $\hat{\zeta}_K (s)$ are the so called non-trivial zeroes of $\zeta_K (s)$, i.e. those in the critical strip $0 < \RRe s < 1$. The famous Riemann conjecture for $K$ asserts that they all lie on the line $\RRe s = 1/2$. \\
Apart from its zeroes, the special values of $\zeta_K (s)$, i.e.\ the numbers $\zeta_K (n)$ for integers $n  \ge 2$, have received a great deal of attention. 
There are two sets of conjectures expressing these values in cohomological terms. One by Bloch and Kato \cite{BK} which has been verified for all abelian extension $K / \Q$, and another more recent one by Lichtenbaum \cite{Li}. Together with the theory of $\zeta$-functions of curves over finite fields this suggests that the Dedekind zeta function should be cohomological in nature. The rest of this article will be devoted to a thorough discussion of this hypothesis in a broader context.

A natural generalization of the Riemann zeta function to the context of arithmetic geometry is the Hasse--Weil zeta function $\zeta_{\eX} (s)$ of an algebraic scheme $\eX / \Z$
\[
\zeta_{\eX} (s) = \prod_{x \in |\eX|} (1 - N (x)^{-s} )^{-1} \; , \; \RRe s > \dim \eX
\]
where $|\eX |$ is the set of closed points of $\eX$ and $N (x)$ is the number of elements in the residue field of $x$. For $\eX = \spec \Oh_K$ we recover the Dedekind zeta function of $K$. It is expected that $\zeta_{\eX} (s)$ has a meromorphic continuation to $\C$ and, if $\eX$ is regular and proper over $\spec \Z$, that
\[
\hat{\zeta}_{\eX} (s) = \zeta_{\eX} (s) \zeta_{\eX_{\infty}} (s)
\]
has a simple functional equation with respect to the substitution of $s$ by $\dim \eX - s$. Here $\zeta_{\eX_{\infty}} (s)$ is a certain product of $\Gamma$-factors depending on the Hodge structure on the cohomology of $\eX_{\infty} = \eX \otimes \R$. This functional equation is known if $\eX$ is equicharacteristic, i.e. an $\F_p$-scheme for some $p$, by using the Lefschetz trace formula and Poincar\'e duality for $l$-adic cohomology. \\
The present strategy for approaching $\zeta_{\eX} (s)$ was first systematically formulated by Langlands. He conjectured that every Hasse--Weil zeta function is up to finitely many Euler factors the product of automorphic $L$-functions. One could then apply the theory of these $L$-functions which is quite well developed in important cases although by no means in general. For $\eX$ with generic fibre related to Shimura varieties this Langlands program has been achieved in very interesting examples. Another spectacular instance was Wiles' proof with Taylor of modularity for most elliptic curves over $\Q$.\\
The strategy outlined in section 3 of the present article is completely different and much closer to the cohomological methods in characteristic $p$.\\ 
By the work of Deligne \cite{De1}, it is known that for proper regular $\eX / \F_p$ the zeroes (resp. poles) of $\hat{\zeta}_{\eX} (s) = \zeta_{\eX} (s)$ have real parts equal to $\nu / 2$ for odd (resp. even) integers $0 \le \nu \le 2 \dim \eX$, and one may wonder whether the same is true for the completed Hasse Weil zeta function $\hat{\zeta}_{\eX} (s)$ of an arbitrary proper and regular scheme $\eX / \Z$.\\ 
As for the orders of vanishing at the integers, a conjecture of Soul\'e \cite{So} asserts that for $\eX / \Z$ regular, quasi-projective connected and of dimension $d$, we have the formula
\begin{equation}\label{eq:0}
\ord_{s = d-n} \zeta_{\eX} (s) = \sum^{2n}_{i=0} (-1)^{i+1} \dim H^i_{\Mh} (\eX , \Q (n)) \; .
\end{equation}
Here, $H^i_{\Mh} (\eX, \Q (n))$ is the rational motivic cohomology of $\eX$ which can be defined by the formula
\[
H^i_{\Mh} (\eX , \Q (n)) = \Gr^n_{\gamma} (K_{2n-i} (\eX) \otimes \Q) \; .
\]
The associated graded spaces are taken with respect to the $\gamma$-filtration on algebraic $K$-theory. Unfortunately it is not even known, except in special cases, whether the dimensions on the right hand side are finite.

\section{The conjectural cohomological formalism}

In this section we interpret some of the conjectures about zeta-functions in terms of an as yet speculative infinite dimensional cohomology theory. We also describe a number of consequences of this very rigid formalism that can be proved directly. Among these there is a formula which expresses the Dedekind $\zeta$-function as a zeta-regularized product. After giving the definition of regularized determinants in a simple algebraic setting we first discuss the formalism in the case of the Dedekind zeta function and then generalize to Hasse--Weil zeta functions.

Given a $\C$-vector space $H$ with an endomorphism $\Theta$ such that $H$ is the countable sum of finite dimensional $\Theta$-invariant subspaces $H_{\alpha}$, the spectrum $\ssp (\Theta)$ is defined as the union of the spectra of $\Theta$ on $H_{\alpha}$, the eigenvalues being counted with their algebraic multiplicities. The (zeta-)regularized determinant $\det_{\infty} (\Theta \tei H)$ of $\Theta$ is defined to be zero if $0 \in \ssp (\Theta)$, and by the formula
\begin{equation}
  \label{eq:1}
  \ddet_{\infty} (\Theta \tei H) := \RRp{\alpha \in \ssp (\Theta)}{} \alpha := \exp (-\zeta'_{\Theta} (0))
\end{equation}
if $0 \notin \ssp (\Theta)$. Here
\[
\zeta_{\Theta} (s) = \sum_{0 \neq \alpha \in \ssp (\Theta)} \alpha^{-s} \; , \quad \mbox{where} \quad - \pi < \mathrm{arg}\, \alpha \le \pi \; ,
\]
is the spectral zeta function of $\Theta$. For (\ref{eq:1}) to make sense we require that $\zeta_{\Theta}$ be convergent in some right half plane, with meromorphic continuation to $\RRe s > - \varepsilon$, for some $\varepsilon > 0$, holomorphic at $s = 0$. For an endomorphism $\Theta_0$ on a real vector space $H_0$, such that $\Theta = \Theta_0 \otimes \id$ on $H = H_0 \otimes \C$ satisfies the above requirements, we set
\[
\ddet_{\infty} (\Theta_0 \tei H_0) = \ddet_{\infty} (\Theta \tei H) \; .
\]
On a finite dimensional vector space $H$ we obtain the ordinary determinant of $\Theta$. As an example of a regularized determinant, consider an endomorphism $\Theta$ whose spectrum consists of the number $1 , 2 , 3 , \ldots$ with multiplicities one. Then
\[
\ddet_{\infty} (\Theta \tei H) = \Rpp{\nu = 1}{\infty} \nu = \sqrt{2\pi} \quad \mbox{since} \quad \zeta' (0) = - \log \sqrt{2 \pi} \; .
\]
The regularized determinant plays a role for example in Arakelov theory and in string theory. In our context it allows us to write the different Euler factors  of zeta-functions in a uniform way as we will first explain for the Dedekind zeta function.\\

Let $\Rh_{\ep}$ for $\ep \nmid \infty$ be the $\R$-vector space of real valued finite Fourier series on $\R / (\log N \ep) \Z$ and set
\[
\Rh_{\ep} = \left\{ \begin{array}{ccl} 
\R [ \exp (-2y)] & \mbox{for real} & \ep \tei \infty \\
\R [\exp (-y)] & \mbox{for complex} & \ep \tei \infty
\end{array} \right.
\]
These spaces of functions (in the variable $y$) carry a natural $\R$-action $\sigma^t$ via $(\sigma^t f) (y) = f (y+t)$ with infinitesimal generator $\Theta = d / dy$. The eigenvalues of $\Theta$ on $\Ch_{\ep} = \Rh_{\ep} \otimes \C$ are just the poles of $\zeta_{\ep} (s)$.

\begin{prop} \label{t1} We have
  $\zeta_{\ep} (s) = \det_{\infty} \left( \frac{1}{2 \pi} (s - \Theta) \tei \Rh_{\ep} \right)^{-1}$ for all places $\ep$ of $K$.
\end{prop}

\begin{proof}
Recall the Hurwitz zeta function $\zeta (s,z)$ which is defined for $\RRe s > 1$ and $z \neq 0, -1 , -2, \ldots$ by the series
\[
\zeta (s,z) = \sum^{\infty}_{\nu = 0} \frac{1}{(z+ \nu)^s} \quad \mbox{where} \; - \pi < \arg (z + \nu) \le \pi \; .
\]
It is known that $\zeta (s,z)$ has an analytic continuation to all $s$ in $\C \ohne \{ 1 \}$ with
\[
\zeta (0,z) = \halb - z \quad \mbox{and} \quad \partial_s \zeta (0,z) = \log \frac{1}{\sqrt{2\pi}} \Gamma (z)
\]
for a suitable branch of $\log \frac{1}{\sqrt{2\pi}} \Gamma (z)$. The latter formula is due to Lerch.

For a complex number $\gamma \neq 0$ we introduce the function
\[
\zeta_{\gamma} (s,z) = \sum^{\infty}_{\nu = 0} \frac{1}{(\gamma (z + \nu))^s} \quad \mbox{where} \; - \pi < \arg \gamma (z+\nu) \le \pi \; .
\]
If $\gamma \neq 0$ is not a negative real number, than we have
\[
\Arg (\gamma (z+\nu)) = \Arg \gamma + \Arg (z + \nu) \quad \mbox{for almost all} \; \nu \ge 0
\]
where $\Arg \in (-\pi , \pi]$ is the prinzipal branch of the argument. Note here that
\[
\lim_{\nu \to \infty} \Arg (z + \nu) = 0 \quad \mbox{and} \quad - \pi < \Arg \gamma < \pi \; .
\]
Hence we have
\[
\zeta_{\gamma} (s,z) = \gamma^{-s} \tzeta (s,z)
\]
where $\tzeta (s,z)$ differs from the Hurwitz zeta function in that finitely many argument in its definition are possibly nonprincipal. Therefore, we still have
\[
\tzeta (0,z) = \halb - z \quad \mbox{and} \quad \exp (-\partial_s \tzeta (0,z)) = \left( \frac{1}{\sqrt{2 \pi}} \Gamma (z) \right)^{-1}
\]
and hence
\[
\zeta_{\gamma} (0,z) = \halb - z \quad \mbox{and} \quad \exp (-\partial_s \zeta_{\gamma} (0,z)) = \gamma^{1/2 -z} \left( \frac{1}{\sqrt{2 \pi}} \Gamma (z) \right)^{-1}
\]
i.e. 
\[
\Rpp{\nu = 0}{\infty} \gamma (z+\nu) = \gamma^{1/2 -z} \left( \frac{1}{\sqrt{2\pi}} \Gamma (z) \right)^{-1} \; .
\]
In order to calculate $\RRP{\nu \in \Z}{} \gamma (z + \nu)$ note that
\[
\sum_{\nu \in \Z} \frac{1}{(\gamma (z + \nu))^s} = \zeta_{\gamma} (s,z) + \zeta_{-\gamma} (s,-z) - (\gamma z)^{-s} \; .
\]
Using the formula
\[
\frac{1}{z} \left( \frac{1}{\sqrt{2 \pi}} \Gamma (z) \right)^{-1} \left( \frac{1}{\sqrt{2\pi}} \Gamma (-z) \right)^{-1} = i (e^{i \pi z} - e^{-i\pi z})
\]
which follows from the equation:
\[
\Gamma (z) \Gamma (1-z) = \frac{\pi}{\sin \pi z}
\]
we get:
\[
\Rpp{\nu \in \Z}{} \gamma (z+\nu) = \left\{ 
\begin{array}{cll}
1 - e^{-2\pi iz} & \mbox{if} & \Imm \gamma > 0 \\
1 - e^{2 \pi iz} & \mbox{if} & \Imm \gamma < 0 
\end{array} \right\} \; .
\]
For $\ep \nmid \infty$ we have:
\begin{eqnarray*}
  \ddet_{\infty} \left( \frac{1}{2 \pi} (s - \theta) \tei \Rh_{\ep} \right) & = & \Rpp{\nu \in \Z}{} \frac{1}{2\pi} \left( s - \frac{2\pi i \nu}{\log N\ep} \right) \\
& = & \Rpp{\nu \in \Z}{} \frac{i}{\log N\ep} \left( \frac{s \log N\ep}{2 \pi i} - \nu \right) \\
& = & 1 - \exp \left( - 2\pi i \frac{s \log N\ep}{2 \pi i} \right) \\
& = & 1 - N\ep^{-s} = \zeta_{\ep} (s)^{-1}
\end{eqnarray*}
for real $\ep \tei \infty$ on the other hand we find
\begin{eqnarray*}
  \ddet_{\infty} \left( \frac{1}{2 \pi} (s - \theta) \tei \Rh_{\ep} \right) & = & \Rpp{\nu = 0}{\infty} \frac{1}{2 \pi} (s + 2 \nu) \\
& = & \Rpp{\nu = 0}{\infty} \frac{1}{\pi} \left(\frac{s}{2} + \nu \right) \\
& = & \pi^{\frac{s}{2} - \halb} \left( \frac{1}{\sqrt{2\pi}} \Gamma \left( \frac{s}{2} \right) \right)^{-1} \\
& = & \sqrt{2} \pi^{s/2} \Gamma \left( \frac{s}{2} \right)^{-1} = \zeta_{\ep} (s)^{-1} \; .
\end{eqnarray*}

Similarly, for complex $\ep \tei \infty$ we get
\begin{eqnarray*}
  \ddet_{\infty} \left( \frac{1}{2 \pi i} (s- \theta) \tei \Rh_{\ep} \right) & = & \Rpp{\nu = 0}{\infty} \frac{1}{2 \pi} (s + \nu) \\
 & = & (2 \pi)^{s - \halb} \left( \frac{1}{\sqrt{2 \pi}} \Gamma (s) \right)^{-1} \\
 & = & (2 \pi)^s \Gamma (s)^{-1} = \zeta_{\ep} (s)^{-1} \; .
\end{eqnarray*}
\end{proof}

In a sense $\overline{\spec \eo} = \spec \eo \cup \{ \ep \tei \infty \}$ is analogous to a projective curve $\eX$ over a finite field of characteristic $p$. For every constructible $\Q_l$-sheaf $\Rh$ on $\eX$ with $l \neq p$ one defines the zeta function of $\eX$ and $\Rh$ by the Euler product over the closed points $x$ of $\eX$
\[
\hzeta_{\eX} (s, \Rh) = \prod_x \ddet (1 - N (x)^{-s} \Fr^*_{\ox} \tei \Rh_{\ox})^{-1} \quad \mbox{for} \; \RRe s > 1 \; .
\]
Here, $\Rh_{\ox}$ is the stalk of $\Rh$ in a geometric point $\ox$ over $x$ and $\Fr_{\ox}$ is the Frobenius morphism. Using the Grothendieck--Verdier--Lefschetz trace formula in characteristic $p$ one gets the following cohomological expression for $\hzeta_{\eX} (s, \Rh)$:
\[
\hzeta_{\eX} (s, \Rh) = \prod^2_{i=0} \ddet (1 - p^{-s} \Fr^*_p \tei H^i_{\et} (\oeX , \Rh))^{(-1)^{i+1}} \; .
\]
Here, $\Fr_p$ is the absolute Frobenius morphism and $\oeX = \eX \otimes_{\F_p} \OF_p$. This formula together with the proposition suggests by analogy that a formula of the following type might hold:
\begin{equation}
  \label{eq:2}
  \hat{\zeta}_K (s) = \prod^2_{i=0} \ddet_{\infty} \textstyle \left( \frac{1}{2 \pi} (s -\Theta) \tei H^i_{\dyn} ( \overline{\spec \eo}  , \Rh) \right)^{(-1)^{i+1}} \; .
\end{equation}
Here $H^i_{\dyn} ( \overline{\spec \eo}  , \Rh) = H^i ((\overline{\spec \eo})_{\dyn} , \Rh)$ would be the cohomology of a space or Grothendieck site $(\overline{\spec \eo})_{\dyn} $
corresponding to $\overline{\spec \eo}$, with coefficients in a sheaf $\Rh$ of $\R$-algebras on $(\overline{\spec \eo})_{\dyn}$. This cohomology should be equipped with a canonical endomorphism $\Theta$.

As recalled earlier $\hat{\zeta}_K (s)$ has poles only at $s = 0,1$ and these are of first order. Moreover the zeroes of $\hat{\zeta}_K (s)$ are just the
non-trivial zeroes of $\zeta_K (s)$. If we assume that the eigenvalues
of $\Theta$ on $H^i_{\dyn} ( \overline{\spec \eo}  , \Rh)$ are distinct
for $i = 0 , 1 , 2$ it follows therefore that
\begin{itemize}
\item $H^0_{\dyn} ( \overline{\spec \eo}  , \Rh) = \R$ with trivial action of $\Theta$, i.e. $\Theta = 0$,
\item $H^1_{\dyn} ( \overline{\spec \eo}  , \Rh)$ is infinite dimensional, the spectrum of $\Theta$ consisting of the non-trivial zeroes $\rho$ of $\zeta_K (s)$ with their multiplicities,
\item $H^2_{\dyn} ( \overline{\spec \eo}  , \Rh) \cong \R$ but with $\Theta = \id$.
\item For $i > 2$ the cohomologies $H^i_{\dyn} ( \overline{\spec \eo}  , \Rh)$ should vanish.
\end{itemize}
Formula (\ref{eq:2}) implies that
\[
\xi_K (s) := \frac{s}{2 \pi} \frac{(s-1)}{2 \pi} \hat{\zeta}_K (s) = \Rpp{\rho}{} \frac{1}{2\pi} (s-\rho) \; .
\]
This formula turned out to be true \cite{D2}, \cite{SchS}. Earlier a related formula had been observed in \cite{K}. \\
If $H$ is some space with an endomorphism $\Theta$ let us write $H (\alpha)$ for $H$ equipped with the twisted endomorphism $\Theta_{H (\alpha)} = \Theta - \alpha \, \id$. With this notation we expect a canonical ``trace''-isomorphism:
\[
\tr : H^2_{\dyn} ( \overline{\spec \eo}  , \Rh) \silo \R (-1) \; .
\]
In our setting the cup product pairing
\[
\cup : H^i_{\dyn} ( \overline{\spec \eo}  , \Rh) \times H^{2-i}_{\dyn} ( \overline{\spec \eo}  , \Rh) \longrightarrow H^2_{\dyn} ( \overline{\spec \eo}  , \Rh) \cong \R (-1)
\]
induces a pairing for every $\alpha$ in $\C$:
\[
{\small \cup : H^i_{\dyn} ( \overline{\spec \eo}  , \Ch)^{\Theta \sim \alpha} \times H^{2-i}_{\dyn} ( \overline{\spec \eo}  , \Ch)^{\Theta \sim 1-\alpha} \longrightarrow H^2_{\dyn} ( \overline{\spec \eo}  , \Ch)^{\Theta \sim 1} \cong \C \; .}
\]
Here $\Ch$ is the complexification of $\Rh$ and $\Theta \sim \alpha$ denotes the subspace of
\[
H^i_{\dyn} ( \overline{\spec \eo}  , \Ch) = H^i_{\dyn} ( \overline{\spec \eo}  , \Rh) \otimes \C
\]
of elements annihilated by some power of $\Theta - \alpha$. We expect Poincar\'e duality in the sense that these pairings should be non-degenerate for all $\alpha$. This is compatible with the functional equation of $\hat{\zeta}_K (s)$. For the precise relation see \cite{D6} 7.19.\\
In the next section we will have more to say on the type of cohomology theory that might be expected for $H^i_{\dyn}  ( \overline{\spec \eo}  , \Rh)$. But first let us note a nice consequence our approach would have. Consider the linear flow $\lambda^t = \exp t \Theta$ on $H^i_{\dyn} ( \overline{\spec \eo}  , \Rh)$. It is natural to expect that it is the flow induced on cohomology by a flow $\phi^t$ on the underlying space $ (\overline{\spec \eo})_{\dyn}$, i.e. $\lambda^t = (\phi^t)^{\ast}$. This implies that $\lambda^t$ would respect cup product and that $\Theta$ would behave as a derivation. Now assume that as in the case of compact Riemann surfaces there is a Hodge $\ast$-operator:
\[
\ast : H^1_{\dyn} ( \overline{\spec \eo}  , \Rh) \silo H^1_{\dyn} ( \overline{\spec \eo}  , \Rh) \; ,
\]
such that
\[
(f , f') = \tr (f \cup (\ast f')) \quad \mbox{for} \; f , f' \; \mbox{in} \; H^1_{\dyn} ( \overline{\spec \eo}  , \Rh) \; ,
\]
is positive definite, i.e. a scalar product on $H^1_{\dyn} ( \overline{\spec \eo}  , \Rh)$. It is natural to assume that $(\phi^t)^{\ast}$ and hence $\Theta$ commutes with $\ast$ on $H^1_{\dyn} ( \overline{\spec \eo}  , \Rh)$. From the equality:
\[
f_1 \cup f_2 = \Theta (f_1 \cup f_2) = \Theta f_1 \cup f_2 + f_1 \cup \Theta f_2
\]
for $f_1 , f_2$ in $H^1_{\dyn} ( \overline{\spec \eo}  , \Rh)$ we would thus obtain the formula
\[
(f_1 , f_2) = (\Theta f_1 , f_2) + (f_1 , \Theta f_2) \; ,
\]
and hence that $\Theta = \halb + A$ where $A$ is a skew-symmetric endomorphism of $H^1_{\dyn} ( \overline{\spec \eo}  , \Rh)$. Hence the Riemann conjecture for $\zeta_K (s)$ would follow. \\
The formula $\Theta = \halb + A$ is also in accordance with numerical investigations on the fluctuations of the spacings between consecutive non-trivial zeroes of $\zeta (s)$. It was found that their statistics resembles that of the fluctuations in the spacings of consecutive eigenvalues of random real skew symmetric matrices, as opposed to the different statistics for random real symmetric matrices; see \cite{KS} for a full account of this story. In fact the comparison was made between hermitian and symmetric matrices, but as pointed out to me by M. Kontsevich, the statistics in the hermitian and real skew symmetric cases agree. \\
The completion of $H^1_{\dyn} ( \overline{\spec \eo}  , \Rh)$ with respect to $( , )$, together with the unbounded operator $\Theta$ would be the space that Hilbert was looking for, and that Berry \cite{B} suggested to realize in a quantum physical setting.

The following considerations are necessary for comparison with the dynamical picture.\\
Formula (\ref{eq:3}) is closely related to a reformulation of the
explicit formulas in analytic number theory using the conjectural
cohomology theory above, see \cite{I} Kap. 3 and \cite{JL} for the
precise relationship. 
\begin{prop}
  \label{t2}
For a test function $\varphi \in \Dh (\R) = C^{\infty}_0 (\R)$ define an entire function $\Phi (s)$ by the formula
\[
\Phi (s) = \int_{\R} \varphi (t) e^{ts} \, dt \; .
\]
Then we have the ``explicit formula'':

\begin{eqnarray}
  \label{eq:3}
  \lefteqn{\Phi (0) - \sum_{\rho} \Phi (\rho) + \Phi (1) = - \log |d_{K / \Q}| \varphi (0)} \nonumber \\
& & + \sum_{\ep \nmid \infty} \log N \ep \left( \sum_{k \ge 1} \varphi (k \log N \ep) + \sum_{k \le -1} N \ep^k \varphi (k \log N \ep) \right) \nonumber \\
& & + \sum_{\ep \tei \infty} W_{\ep} (\varphi) \; .
\end{eqnarray}
\end{prop}

Here $\rho$ runs over the non-trivial zeroes of $\zeta_K (s)$ i.e. those that are contained in the critical strip $0 < \RRe s < 1$. Moreover $\ep$ runs over the places of $K$ and $d_{K / \Q}$ is the discriminant of $K$ over $\Q$. For $\ep \tei \infty$ the $W_{\ep}$ are distributions which are determined by the $\Gamma$-factor at $\ep$. If $\varphi$ has support in $\R^{> 0}$ then
\[
W_{\ep} (\varphi) = \int^{\infty}_{-\infty} \frac{\varphi (t)}{1 - e^{\kappa_{\ep} t}} \, dt
\]
where $\kappa_{\ep} = -1$ if $\ep$ is complex and $\kappa_{\ep} = -2$ if $\ep$ is real. If $\varphi$ has support on $\R^{< 0}$ then
\[
W_{\ep} (\varphi) = \int^{\infty}_{-\infty} \frac{\varphi (t)}{1 - e^{\kappa_{\ep} |t|}} \, e^t \, dt \; .
\]
There are different ways to write $W_{\ep}$ on all of $\R$ but we will not discuss this here. See for example \cite{Ba} which also contains a proof of the theorem for much more general test functions.

Formula (\ref{eq:3}) can be written equivalently as an equality of distributions on $\R$:
\begin{eqnarray}
  \label{eq:4}
  \lefteqn{1 - \sum_{\rho} e^{t\rho} + e^t  =  - \log |d_{K / \Q}| \delta_0} \nonumber \\
 & &  + \sum_{\ep \nmid \infty} \log N\ep \left( \sum_{k \ge 1} \delta_{k \log N \ep} + \sum_{k \le -1} N \ep^{k} \delta_{k \log N \ep} \right) \nonumber \\
& & + \sum_{\ep \tei \infty} W_{\ep} \; .
\end{eqnarray}

For fixed $t$, the numbers $e^{t\rho}$ should be the eigenvalues of $\phi^{t*} = \exp t\Theta$ on $H^1_{\dyn} (\overline{\spec \eo} , \Ch)$. Hence $\sum_{\rho} e^{t\rho}$ should be the trace of $\phi^{t*}$ on this cohomology. However this series does not converge for any real $t$ since the numbers $e^{t\rho}$ do not tend to zero. However, as a distribution on $\R$ the series $\sum_{\rho} e^{t\rho}$ converges since the series $\sum_{\rho} \langle e^{t\rho} , \varphi \rangle = \sum_{\rho} \Phi (\rho)$ converges for every test function $\varphi$. We therefore define the distributional trace of $\phi^*$ on $H^1_{\dyn} (\overline{\spec \eo_K} , \Rh)$ to be the distribution $\sum_{\rho} e^{t\rho}$ on $\R$. Similarly $1$ resp. $e^t$ will be the distributional trace of $\phi^*$ on $H^i_{\dyn} (\overline{\spec \eo} , \Rh)$ for $i = 0$ resp. $i = 2$. 
 Conjecturally (\ref{eq:4}) can thus be reformulated as the following identity of distributions
\begin{eqnarray}
  \label{eq:5}
  \lefteqn{\sum_i (-1)^i \Tr (\phi^{\ast} \tei H^i_{\dyn} ( \overline{\spec \eo}  , \Rh))  =  - \log |d_{K / \Q}| \delta_0 } \nonumber \\
& & + \sum_{\ep \nmid \infty} \log N \ep \left( \sum_{k \ge 1} \delta_{k \log N\ep} + \sum_{k \le -1} N \ep^k \delta_{k \log N\ep} \right) \nonumber \\
& & + \sum_{\ep \tei \infty} W_{\ep} \; .
\end{eqnarray}

We now turn to Hasse--Weil zeta functions of algebraic schemes $\eX / \Z$. A similar argument as for the Riemann zeta function suggests that with $d = \dim \eX$ we have
\begin{equation} \label{eq:6}
\zeta_{\eX} (s) = \prod^{2d}_{i=0} \ddet_{\infty} \left( \frac{1}{2 \pi} (s-\Theta) \tei H^i_{\dyn,c} ( \eX  , \Rh) \right)^{(-1)^{i+1}}
\end{equation}
where $H^i_{\dyn, c} ( \eX  , \Rh)$ is some real cohomology with compact supports associated to a dynamical system $\eX_{\dyn}$, i.e. a space or a site with an $\R$-action $\phi^t$ attached to $\eX$. Here $\Theta$ should be the infinitesimal generator of the induced flow $\phi^{t*}$ on cohomology. In particular we would have
\[
\ord_{s = \alpha} \; \zeta_{\eX} (s) \; = \; \sum^{2d}_{i=0} (-1)^{i+1} \dim H^i_{\dyn, c} ( \eX  , \Ch)^{\Theta \sim \alpha} \; .
\]
Now let $\eX$ be a regular scheme which is projective over $\spec \Z$ and equidimensional of dimension $d$. Then Poincar\'e duality
\begin{equation} \label{eq:7}
\cup : H^i_{\dyn, c} ( \eX  , \Rh) \times H^{2d-i}_{\dyn} ( \eX  , \Rh) \longrightarrow H^{2d}_{\dyn,c} ( \eX  , \Rh) \silo \R (-d)
\end{equation}
should identify 
\[
H^i_{\dyn, c} ( \eX  , \Ch)^{\Theta \sim \alpha} \quad \mbox{with the dual of} \quad H^{2d-i}_{\dyn} ( \eX  , \Ch)^{\Theta \sim d - \alpha} \; .
\]
In particular we would get:
\[
\ord_{s = d-n} \; \zeta_{\eX} (s) \; = \; \sum^{2d}_{i=0} (-1)^{i+1} \dim H^i_{\dyn} ( \eX  , \Ch (n))^{\Theta \sim 0} \; ,
\]
where $\Ch (\alpha)$ is the sheaf $\Ch$ on $ \eX $ with action of the flow twisted by $e^{-\alpha t}$. Thus
\[
H^i_{\dyn} ( \eX  , \Ch (n))^{\Theta \sim 0} = H^i_{\dyn} ( \eX  , \Ch)^{\Theta \sim n} \; .
\]
We also expect formal analogues of Tate's conjecture
\begin{equation}
  \label{eq:8}
  H^i_{\Mh} (\eX , \C (n)) := \Gr^n_{\gamma} K_{2n-i} (\eX) \otimes \C \silo H^i_{\dyn} ( \eX  , \Ch (n))^{\Theta \sim 0} \; ,
\end{equation}
and in particular that
\[
H^i_{\dyn} ( \eX  , \Ch (n))^{\Theta\sim 0} = 0 \quad \mbox{for} \; i > 2n \; .
\]
These assertions would imply Soul\'e's conjecture (\ref{eq:0}). \\
In terms of cohomology, the explicit formulas would take the form
\begin{eqnarray}
  \label{eq:9}
  \lefteqn{\sum_i (-1)^i \Tr (\phi^{\ast} \tei H^i_{\dyn, c} ( \eX  , \Rh))_{\diss} } \\ 
  & = & - (\log A_{\eX}) \delta_0 + \sum_{x \in |\eX|} \log N (x) \sum_{k\ge 1} \delta_{k \log N (x)} \nonumber \\
 & & + \sum_{x \in |\eX|} \log N (x) \sum_{k \le -1} N (x)^k \delta_{k \log N (x)} \nonumber  \; .
\end{eqnarray}
Here $A_{\eX}$ is the conductor of $\eX$.

In support of these ideas we have the following result.
\begin{theorem} \label{t3} On the category of algebraic $\F_p$-schemes $\eX$ there is a cohomology theory in $\C$-vector spaces with a linear flow such that {\rm (\ref{eq:6})} holds. For a regular connected projective $\eX$ of equidimension $d$ it satisfies Poincar\'e duality {\rm (\ref{eq:7})}. Moreover {\rm (\ref{eq:8})} reduces to the Tate conjecture for $l$-adic cohomology.
\end{theorem}
See \cite{D3} \S\,2, \cite{D6} \S\,4 for more precise statements and the simple construction based on $l$-adic cohomology. That construction cannot be generalized to schemes $\eX$ which are flat over $\spec \Z$.

If there were a dynamical cohomology theory $H^i_{\dyn} ( \overline{\eX}  , \Rh)$ attached to some Arakelov compactification $\overline{\eX}$ of $\eX$ such that
\[
\hat{\zeta}_{\eX} (s) = \prod^{2d}_{i=0} \ddet_{\infty} \left( \frac{1}{2\pi} (s - \Theta) \tei H^i_{\dyn} ( \overline{\eX}  , \Rh) \right)^{(-1)^{i+1}} \; ,
\]
then as above Poincar\'e duality for $H^i_{\dyn} ( \overline{\eX}  , \Rh)$ would be in accordance with the expected functional equation for $\hat{\zeta}_{\eX} (s)$. A Hodge $\ast$-operator
\[
\ast : H^i_{\dyn} ( \overline{\eX}  , \Rh) \longrightarrow H^{2d-i}_{\dyn} ( \overline{\eX}  , \Rh)
\]
defining a scalar product via $(f,f') = \tr (f \cup (\ast f'))$ and for which
\begin{equation}
\label{eq:11n}
  \phi^{t\ast} \verk \ast = (e^t)^{d-i} \ast \verk \phi^{t\ast} \; , \quad \mbox{i.e.} \quad \Theta \verk \ast = \ast \verk (d-i + \Theta) \; ,
\end{equation}
holds, would imply that $\Theta - i/2$ is skew symmetric, hence the Riemann hypotheses for $\hat{\zeta}_{\eX} (s)$. On the level of forms equation (\ref{eq:11n}) means that the flow changes the metric defining the $\ast$-operator by the conformal factor $e^t$. However, equation (\ref{eq:11n}) for forms is a much stronger condition than for cohomology classes.

In \cite{D2} we constructed cohomology $\R$-vector spaces $H^i_{\ar}$ with a linear flow on the category of varieties over $\R$ or $\C$ such that
\[
\zeta_{\eX_{\infty}} (s) = \prod^{2 \dim \eX_{\infty}}_{i=0} \ddet_{\infty} \left( \frac{1}{2 \pi} (s- \Theta) \tei H^i_{\ar} (\eX_{\infty}) \right)^{(-1)^{i+1}} \; .
\]
Cup product and functoriality turn the spaces $H^i_{\ar} (\eX_{\infty})$ into modules under $H^0_{\ar} (\eX_{\infty}) = H^0_{\ar} (\spec \R) = \Rh_{\infty}$ of rank equal to $\dim H^i (\eX_{\infty} , \Q)$. Philosophically the scheme $\eX$ should have bad semistable ``reduction'' at infinity. In accordance with this idea Consani \cite{Cons} has refined the theory $H^i_{\ar}$ to a cohomology theory with a linear flow and a monodromy operator $N$ which contains $H^i_{\ar}$ as the kernel of $N$.

\section{Foliations and their cohomology}
A $d$-dimensional foliation $\Fh = \Fh_X$ of a smooth manifold $X$ of dimension $a$ is a partition of $X$ into immersed connected $d$-dimensional manifolds $F$, the ``leaves''. Locally the induced partition should be trivial: Every point of $X$ should have an open neighborhood $U$ diffeomorphic to an open ball $B$ in $\R^a$ such that the leaves of the induced partition on $U$ correspond to the submanifolds $B \cap (\R^d \times \{ y \})$ of $B$ for $y$ in $\R^{a-d}$.

One of the simplest non-trivial examples is the one-dimensional foliation of the two-dimensional torus $T^2 = \R^2 / \Z^2$ by lines of irrational slope $\alpha$. These are given by the immersions
\[
\R \hookrightarrow T^2 \; , \; t \mapsto (x + t \alpha , t) \mod \Z^2
\]
parametrized by $x \mod \Z + \alpha \Z$. In this case every leaf is dense in $T^2$ and the intersection of a global leaf with a small open neighborhood $U$ as above decomposes into countably many connected components. It is the global behaviour which makes foliations complicated. For an introduction to foliation theory, the reader may turn to \cite{Go} for example.

Let $\Rh = \Rh_{\Fh}$ be the sheaf of germs of smooth real valued functions which are locally constant on the leaves. We are interested in the cohomology groups $H^i (X, \Rh)$. 

To a foliation $\Fh$ on $X$ we may attach its tangent bundle $T \Fh$ whose total space is the union of the tangent spaces to the leaves. By local triviality of the foliation it is a sub vector bundle of the tangent bundle $TX$. It is integrable i.e. the commutator of any two vector fields with values in $T \Fh$ again takes values in $T\Fh$. Conversely a theorem of Frobenius asserts that every integrable sub vector bundle of $TX$ arises in this way.

On an open subset $U$ of $X$ the differential forms of order $n$ along the leaves are defined as the smooth sections of the real vector bundle $\Lambda^n T^* \Fh$, 
\[
\Ah^n_{\Fh} (U) = \Gamma (U, \Lambda^n T^* \Fh) \; .
\]
The same formulas as in the classical case define exterior derivatives along the leaves:
\[
d^n_{\Fh} : \Ah^n_{\Fh} (U) \longrightarrow \Ah^{n+1}_{\Fh} (U)
\]
which satisfy the relation $d^{n+1}_{\Fh} \verk d^n_{\Fh} = 0$. 

The Poincar\'e lemma extends to the foliated context and one obtains a resolution of $\Rh$ by the fine sheaves $\Ah^n_{\Fh}$:
\[
0 \longrightarrow \Rh \longrightarrow \Ah^0_{\Fh} \xrightarrow{d_{\Fh}} \Ah^1_{\Fh} \xrightarrow{d_{\Fh}} \ldots
\]
Hence we have the following de Rham description of the cohomology of $\Rh$:
\[
H^n (X, \Rh) = \Ker (d^n_{\Fh} : \Ah^n_{\Fh} (X) \to \Ah^{n+1}_{\Fh} (X)) / \Imm (d^{n-1}_{\Fh} : \Ah^{n-1}_{\Fh} (X) \to \Ah^n_{\Fh} (X)) \; .
\]
For our purposes these invariants are actually too subtle. We therefore consider the reduced leafwise cohomology
\[
\oH^n (X , \Rh) = \Ker d^n_{\Fh} / \overline{\Imm d^{n-1}_{\Fh}} \; .
\]
Here the quotient is taken with respect to the topological closure of $\Imm d^{n-1}_{\Fh}$ in the natural Fr\'echet topology on $\Ah^n_{\Fh} (X)$. The reduced cohomologies are nuclear Fr\'echet spaces. Even if the leaves are dense, already $\oH^1 (X,\Rh)$ can be infinite dimensional, c.f. \cite{DS1}.

The cup product pairing induced by the exterior product of forms along the leaves turns $\oH^{\hullet} (X,\Rh)$ into a graded commutative $H^0 (X,\Rh)$-algebra.

For the torus foliation above with $\alpha \notin \Q$ we have $\oH^0 (T^2 , \Rh) = \R$. Some Fourier analysis reveals that $\oH^1 (T^2 , \Rh) \cong \R$. The higher cohomologies vanish since by the de\,Rham description we have
\[
H^n (X, \Rh) = 0 \quad \mbox{for all} \; n > d = \dim \Fh \; .
\]
For a smooth map $f : X \to Y$ of foliated manifolds which maps leaves into leaves, continuous pullback maps
\[
f^* : \Ah^n_{\Fh_Y} (Y) \longrightarrow \Ah^n_{\Fh_X} (X)
\]
are defined for all $n$. They commute with $d_{\Fh}$ and respect the exterior product of forms. Hence they induce a continuous map of reduced cohomology algebras
\[
f^* : \oH^{\hullet} (Y , \Rh_Y) \longrightarrow \oH^{\hullet} (X , \Rh_X) \; .
\]
A (complete) flow is a smooth $\R$-action $\phi : \R \times X \to X , (t,x) \mapsto \phi^t (x)$. It is called $\Fh$-compatible if every diffeomorphism $\phi^t : X \to X$ maps leaves into leaves. If this is the case we obtain a linear $\R$-action $t \mapsto \phi^{t*}$ on $\oH^n (X , \Rh)$ for every $n$. Let
\[
\Theta : \oH^n (X, \Rh) \longrightarrow \oH^n (X,\Rh)
\]
denote the infinitesimal generator of $\phi^{t*}$:
\[
\Theta h = \lim_{t\to 0} \frac{1}{t} (\phi^{t*} h -h ) \; .
\]
The limit exists and $\Theta$ is continuous in the Fr\'echet topology. As $\phi^{t*}$ is an algebra endomorphism of the $\R$-algebra $\oH^{\hullet} (X,\Rh)$ it follows that $\Theta$ is an $\R$-linear derivation. Thus we have
\begin{equation}
  \label{eq:10}
  \Theta (h_1 \cup h_2) = \Theta h_1 \cup h_2 + h_1 \cup \Theta h_2
\end{equation}
for all $h_1 , h_2$ in $\oH^n (X,\Rh)$.

For arbitrary foliations the reduced leafwise cohomology does not seem to have a good structure theory. For Riemannian foliations however the situation is much better. These foliations are characterized by the existence of a ``bundle-like'' metric $g$. This is a Riemannian metric whose geodesics are perpendicular to all leaves whenever they are perpendicular to one leaf. For example any one-codimensional foliation given by a closed one-form without singularities is Riemannian. 

Assuming that $X$ is also oriented, the graded Fr\'echet space $\Ah^{\hullet}_{\Fh} (X)$ carries a canonical inner product:
\[
(\alpha , \beta) = \int_X \langle \alpha , \beta \rangle_{\Fh} \vol \; .
\]
Here $\langle , \rangle_{\Fh}$ is the Riemannian metric on $\Lambda^{\hullet} T^* \Fh$ induced by $g$ and $\vol$ is the volume form or density on $X$ coming from $g$. Let
\[
\Delta_{\Fh} = d_{\Fh} d^*_{\Fh} + d^*_{\Fh} d_{\Fh}
\]
denote the Laplacian using the formal adjoint of $d_{\Fh}$ on $X$. Since $\Fh$ is Riemannian the restriction of $\Delta_{\Fh}$ to any leaf $F$ is the Laplacian on $F$ with respect to the induced metric \cite{AK1} Lemma 3.2, i.e.
\[
(\Delta_{\Fh} \alpha) \, |_F = \Delta_F (\alpha \, |_F) \quad \mbox{for all} \; \alpha \in \Ah^{\hullet}_{\Fh} (X) \; .
\]
We now assume that also $T \Fh$ is orientable. Via $g$ the choice of an orientation determines a volume form $\vol_{\Fh}$ in $\Ah^d_{\Fh} (X)$ and hence a Hodge $*$-operator
\[
*_{\Fh} : \Lambda^n T^*_x \Fh \silo \Lambda^{d-n} T^*_x \Fh \quad \mbox{for every} \; x \; \mbox{in} \; X \; .
\]
It is determined by the condition that
\[
v \wedge *_{\Fh} w = \langle v,w \rangle_{\Fh} \, \vol_{\Fh , x} \quad \mbox{for} \; v,w \; \mbox{in} \; \Lambda^{\hullet} T^*_x \Fh \; .
\]
These fibrewise star-operators induce the leafwise $*$-operator on forms:
\[
*_{\Fh} : \Ah^n_{\Fh} (X) \silo \Ah^{d-n}_{\Fh} (X) \; .
\]
We now list some important properties of leafwise cohomology.

{\bf Properties} Assume that $X$ is compact, $\Fh$ a $d$-dimensional oriented Riemannian foliation and $g$ a bundle-like metric for $\Fh$.

Then the natural map
\begin{equation}
  \label{eq:11}
  \Ker \Delta^n_{\Fh} \silo \oH^n (X , \Rh ) \; ,\; \omega \longmapsto \omega \bmod \overline{\Imm d^{n-1}_{\Fh}}
\end{equation}
is a topological isomorphism of Fr\'echet spaces. We denote its inverse by $\Hh$.

This result is due to \'Alvarez L\'opez and Kordyukov \cite{AK1}. It is quite deep since $\Delta_{\Fh}$ is only elliptic along the leaves so that the ordinary elliptic regularity theory does not suffice. For non-Riemannian foliations (\ref{eq:11}) does not hold in general \cite{DS1}. All the following results are consquences of this Hodge theorem.

The Hodge $*$-operator induces an isomorphism 
\[
*_{\Fh} : \Ker \Delta^n_{\Fh} \silo \Ker \Delta^{d-n}_{\Fh}
\]
since it commutes with $\Delta_{\Fh}$ up to sign. From (\ref{eq:11}) we therefore get isomorphisms for all $n$:
\begin{equation}
  \label{eq:12}
  *_{\Fh} : \oH^n (X,\Rh) \silo \oH^{d-n} (X,\Rh) \; .
\end{equation}
For the next property define the trace map
\[
\tr : \oH^d (X,\Rh) \longrightarrow \R
\]
by the formula
\[
\tr (h) = \int_X *_{\Fh} (\Hh (h)) \vol \; .
\]
It is an isomorphism if $\Fh$ has a dense leaf. Note that for {\it any} representative $\alpha$ in the cohomology class $h$ we have
\[
\tr (h) = \int_X *_{\Fh} (\alpha) \vol \; .
\]
Namely $\alpha - \Hh (h) = d_{\Fh} \beta$ and
\begin{eqnarray*}
  \int_X *_{\Fh} (d_{\Fh} \beta) \vol & = & \pm \int_X d^*_{\Fh} (*_{\Fh} \beta) \vol \\
& = &\pm ( 1 , d^*_{\Fh} (*_{\Fh} \beta)) \\
& = & \pm ( d_{\Fh} (1) , *_{\Fh} \beta)\\
& = & 0 \; .
\end{eqnarray*}
It is not difficult to see using (\ref{eq:11}) that we get a scalar product on $\oH^n (X,\Rh)$ for every $n$ by setting:
\begin{eqnarray}
  \label{eq:13}
  (h,h') & = & \tr (h \cup *_{\Fh} h') \\
& = & \int_X \langle \Hh (h) , \Hh (h') \rangle_{\Fh} \vol \; . \nonumber
\end{eqnarray}
It follows from this that the cup product pairing 
\begin{equation}
  \label{eq:14}
  \cup : \oH^n (X,\Rh) \times \oH^{d-n} (X,\Rh) \longrightarrow \oH^d (X,\Rh) \xrightarrow{tr} \R
\end{equation}
is non-degenerate.

Next we discuss the K\"unneth formula. Assume that $Y$ is another compact manifold with a Riemannian foliation. Then the canonical map
\[
H^n (X,\Rh) \otimes H^m (Y,\Rh) \longrightarrow H^{n+m} (X \times Y, \Rh)
\]
induces a topological isomorphism \cite{M}:
\begin{equation}
  \label{eq:15}
  \oH^n (X,\Rh) \hat{\otimes} \oH^m (Y,\Rh) \silo \oH^{n+m} (X \times Y,\Rh) \; .
\end{equation}
Since the reduced cohomology groups are nuclear Fr\'echet spaces, it does not matter which topological tensor product is chosen in (\ref{eq:15}). The proof of this K\"unneth formula uses (\ref{eq:11}) and the spectral theory of the Laplacian $\Delta_{\Fh}$.

Before we deal with more specific topics let us mention that Hodge--K\"ahler theory can also be generalized to the foliated context. A complex structure on a foliation $\Fh$ is an almost complex structure $J$ on $T \Fh$ such that all restrictions $J \, |_F$ to the leaves are integrable. Then the leaves carry holomorphic structures which vary smoothly in the transverse direction. A foliation $\Fh$ with a complex structure $J$ is called K\"ahler if there is a hermitian metric $h$ on the complex bundle $T_c \Fh = (T \Fh , J)$ such that the K\"ahler form along the leaves
\[
\omega_{\Fh} = - \halb \imm h \in \Ah^2_{\Fh} (X)
\]
is closed. Note that for example any foliation by orientable surfaces can be given a K\"ahlerian structure by choosing a metric on $X$, c.f. \cite{MS} Lemma A.3.1. Let
\[
L_{\Fh} : \oH^n (X,\Rh) \longrightarrow \oH^{n+2} (X,\Rh) \; , \; L_{\Fh} (h) = h \cup [\omega_{\Fh}]
\]
denote the Lefschetz operator.

The following assertions are consequences of (\ref{eq:11}) combined with the classical Hodge--K\"ahler theory. See \cite{DS3} for details. Let $X$ be a compact orientable manifold and $\Fh$ a K\"ahlerian foliation with respect to the hermitian metric $h$ on $T_c \Fh$. Assume in addition that $\Fh$ is Riemannian. Then we have:
\begin{equation}
  \label{eq:16}
  \oH^n (X,\Rh) \otimes \C = \bigoplus_{p+q=n} H^{pq} \; , \quad \mbox{where} \; \overline{H^{pq}} = H^{qp} \; .
\end{equation}
Here $H^{pq}$ consists of those classes that can be represented by $(p,q)$-forms along the leaves. Moreover there are topological isomorphisms
\[
H^{pq} \cong \oH^q (X , \Omega^p_{\Fh})
\]
with the reduced cohomology of the sheaf of holomorphic $p$-forms along the leaves. 

Furthermore the Lefschetz operator induces isomorphisms
\begin{equation}
  \label{eq:17}
  L^i_{\Fh} : \oH^{d-i} (X,\Rh) \silo \oH^{d+i} (X,\Rh) \quad \mbox{for} \; 0 \le i \le d \; .
\end{equation}
Finally the space of primitive cohomology classes $\oH^n (X,\Rh)_{\prim}$ carries the structure of a polarizable $\ind \, \R$-Hodge structure of weight $n$.

After this review of important properties of the reduced leafwise cohomology of Riemannian foliations we turn to a specific result relating flows and cohomology. This is a kind of toy model for the dynamical cohomologies in the preceeding section. The existence of a conformal metric for the flow simplifies the analysis. However, I do not think that such a metric will exist for dynamical systems relevant to number fields. In order to verify equation (\ref{eq:11n}) for them one will need the K\"ahler identities on cohomology.

\begin{theorem}
  \label{t4}
Let $X$ be a compact $3$-manifold and $\Fh$ a Riemannian foliation by surfaces with a dense leaf. Let $\phi^t$ be an $\Fh$-compatible flow on $X$ which is conformal on $T\Fh$ with respect to a metric $g$ on $T\Fh$ in the sense that for some constant $\alpha$ and all $x \in X$ and $t \in \R$ we have: 
\begin{equation}
  \label{eq:18}
  g (T_x \phi^t (v) , T_x \phi^t (w)) = e^{\alpha t} g (v,w) \; \mbox{for all} \; v,w \in T_x \Fh \; .
\end{equation}
Then we have for the infinitesimal generator of $\phi^{t*}$ that:
\[
\Theta = 0 \; \mbox{on} \; \oH^0 (X,\Rh) = \R \quad \mbox{and} \quad \Theta = \alpha \; \mbox{on} \; \oH^2 (X,\Rh) \cong \R \; .
\]
On $\oH^1 (X,\Rh)$ the operator $\Theta$ has the form
\[
\Theta = \frac{\alpha}{2} + S
\]
where $S$ is skew-symmetric with respect to the inner product $( , )$ above.
\end{theorem}

\begin{remarknn}
   For the bundle-like metric on $X$ required for the construction of $(,)$ we take any extension of the given metric on $T \Fh$ to a bundle-like metric on $TX$. Such extensions exist.
\end{remarknn}

\begin{proofof}
  {\bf \ref{t4}} Because we have a dense leaf, $\oH^0 (X,\Rh) = H^0 (X,\Rh)$ consists only of constant functions. On these $\phi^{t*}$ acts trivially so that $\Theta = 0$. Since
\[
*_{\Fh} : \R = \oH^0 (X,\Rh) \silo \oH^2 (X,\Rh)
\]
is an isomorphism and since
\[
\phi^{t*} (*_{\Fh} (1)) = e^{\alpha t} (*_{\Fh} 1)
\]
by conformality, we have $\Theta = \alpha$ on $\oH^2 (X,\Rh)$.

For $h_1 , h_2$ in $\oH^1 (X,\Rh)$ we find
\begin{equation}
  \label{eq:19}
  \alpha (h_1 \cup h_2) = \Theta (h_1 \cup h_2) = \Theta h_1 \cup h_2 + h_1 \cup \Theta h_2 \; .
\end{equation}
By conformality $\phi^{t*}$ commutes with $*_{\Fh}$ on $\oH^1 (X,\Rh)$. Differentiating, it follows that $\Theta$ commutes with $*_{\Fh}$ as well. Since by definition we have
\[
(h, h') = \tr (h \cup *_{\Fh} h') \quad \mbox{for} \; h , h' \in \oH^1 (X,\Rh) \; ,
\]
it follows from (\ref{eq:19}) that as desired:
\[
\alpha ( h , h' ) = ( \Theta h , h' ) + ( h , \Theta h' ) \; .
\]
\end{proofof}

\section{Dynamical Lefschetz trace formulas}
The formulas we want to consider in this section relate the compact orbits of a flow with the alternating sum of suitable traces on cohomology. A suggestive but non-rigorous argument of Guillemin \cite{G} later rediscovered by Patterson \cite{P} led to the following conjecture \cite{D10} \S\,3. Let $X$ be a compact manifold with a one-codimensional foliation $\Fh$ and an $\Fh$-compatible flow $\phi$. Assume that the fixed points and the periodic orbits of the flow are non-degenerate in the following sense: For any fixed point $x$ and every $t > 0$, the tangent map $T_x \phi^t$ has eigenvalues different from $1$. For any closed orbit $\gamma$ of length $l (\gamma)$ and any $x \in \gamma$ and integer $k \neq 0$ the automorphism $T_x \phi^{kl (\gamma)}$ of $T_x X$ should have the eigenvalue $1$ with algebraic multiplicity one. Observe that the vector field $Y_{\phi}$ generated by the flow provides an eigenvector $Y_{\phi,x}$ for the eigenvalue $1$. 

Recall that the length $l (\gamma) > 0$ of $\gamma$ is defined by the isomorphism:
\[
\R / l (\gamma) \Z \silo \gamma \; , \; t \longmapsto \phi^t (x) \; .
\]
For a fixed point $x$ we set\footnote{This is different from the normalization in \cite{D10} \S\,3.}
\[
\varepsilon_x = \sgn \det (1 - T_x \phi^t \tei T_x \Fh) \; .
\]
This is independent of $t > 0$. For a closed orbit $\gamma$ and $k \in \Z \ohne 0$ set$^1$
\[
\varepsilon_{\gamma} (k) = \sgn \det (1 - T_x \phi^{kl (\gamma)} \tei T_x X / \R Y_{\phi,x}) \; .
\]
It does not depend on the point $x \in \gamma$.

Finally let $\Dh' (J)$ denote the space of Schwartz distributions on an open subset $J$ of $\R$. 

\begin{conj}
  \label{t5}
For $X , \Fh$ and $\phi$ as above there exists a natural definition of a $\Dh' (\R^{>0})$-valued trace of $\phi^*$ on the reduced leafwise cohomology $\oH^{\hullet} (X,\Rh)$ such that in $\Dh' (\R^{>0})$ we have:
\begin{equation}
  \label{eq:20}
\sum\limits^{\dim \Fh}_{n=0} (-1)^n \Tr (\phi^* \tei \oH^n (X,\Rh)) = \sum\limits_{\gamma} l (\gamma) \sum\limits^{\infty}_{k=1} \varepsilon_{\gamma} (k) \delta_{kl (\gamma)} + \sum\limits_x \varepsilon_x |1 - e^{\kappa_x t}|^{-1} \; .
\end{equation}
\end{conj}

Here $\gamma$ runs over the closed orbits of $\phi$ which are not contained in a leaf and $x$ over the fixed points. For $a \in \R , \delta_a$ is the Dirac distribution in $a$ and $\kappa_x$ is defined by the action of $T_x \phi^t$ on the $1$-dimensional vector space $T_x X / T_x \Fh$. That action is multiplication by $e^{\kappa_x t}$ for some $\kappa_x \in \R$ and all $t$.

The conjecture is not known (except for $\dim X = 1$) if $\phi$ has fixed points. It may well have to be amended somewhat in that case. The analytic difficulty in the presence of fixed points lies in the fact that in this case $\Delta_{\Fh}$ cannot be transversally elliptic to the $\R$-action by the flow, so that the methods of transverse index theory do not apply directly. In the simpler case when the flow is everywhere transversal to $\Fh$, \'Alvarez L\'opez and Kordyukov have proved a  beautiful strengthening of the conjecture. Partial results were obtained by other methods in \cite{Laz}, \cite{DS2}. We now describe their result in a convenient way for our purposes:

\begin{punkt}
  \label{t6} \rm
Assume $X$ is a compact oriented manifold with a one codimensional foliation $\Fh$. Let $\phi$ be a flow on $X$ which is everywhere transversal to the leaves of $\Fh$. Then $\Fh$ inherits an orientation and it is Riemannian \cite{Go} III 4.4. Fixing a bundle-like metric $g$, the cohomologies $\oH^n (X,\Rh)$ acquire pre-Hilbert structures (\ref{eq:13}) and we can consider their Hilbert space completions $\hat{H}^n (X,\Rh)$. For every $t$ the linear operator $\phi^{t*}$ is bounded on $(\oH^n (X,\Rh) , \| \; \|)$ and hence can be continued uniquely to a bounded operator on $\hat{H}^n (X,\Rh)$ c.f. theorem \ref{t8}.
\end{punkt}

By transversality the flow has no fixed points. We assume that all periodic orbits are non-degenerate.

\begin{theorem}[\cite{AK2}]
  \label{t7}
Under the conditions of (\ref{t6}), for every test function $\varphi \in \Dh (\R) = C^{\infty}_0 (\R)$ the operator
\[
A_{\varphi} = \int_{\R} \varphi (t) \phi^{t*} \, dt
\]
on $\hat{H}^n (X,\Rh)$ is of trace class. Setting:
\[
\Tr (\phi^* \tei \oH^n (X,\Rh)) (\varphi) = \tr A_{\varphi}
\]
defines a distribution on $\R$. The following formula holds in $\Dh' (\R)$:
\begin{equation}
  \label{eq:21}
  \sum^{\dim \Fh}_{n=0} (-1)^n \Tr (\phi^* \tei \oH^n (X,\Rh)) = \chi_{\Co} (\Fh , \mu) \delta_0 + \sum_{\gamma} l (\gamma) \sum_{k \in \Z \ohne 0} \varepsilon_{\gamma} (k) \delta_{kl (\gamma)} \; .
\end{equation}
Here $\chi_{\Co} (\Fh , \mu)$ denotes Connes' Euler characteristic of the foliation with respect to the transverse measure $\mu$ corresponding to $\tr$. (See \cite{MS}.)
\end{theorem}

It follows from the theorem that if the right hand side of (\ref{eq:21}) is non-zero, at least one of the cohomology groups $\oH^n_{\Fh} (X)$ must be infinite dimensional. Otherwise the alternating sum of traces would be a smooth function and hence have empty singular support.

By the Hodge isomorphism (\ref{eq:11}) one may replace cohomology by the spaces of leafwise harmonic forms. The left hand side of the dynamical Lefschetz trace formula then becomes the $\Dh' (\R)$-valued transverse index of the leafwise de Rham complex. Note that the latter is transversely elliptic for the $\R$-action $\phi^t$. Transverse index theory with respect to compact group actions was initiated in \cite{A}. A definition for non-compact groups of a transverse index was later given by H\"ormander \cite{Si} Appendix II. 

As far as we know the relation of (\ref{eq:21}) with transverse index theory in the sense of Connes--Moscovici still needs to be clarified.

Let us now make some remarks on the operators $\phi^{t*}$ on $\hat{H}^n (X,\Rh)$ in a more general setting:

\begin{theorem}
  \label{t8}
Let $\Fh$ be a Riemannian foliation on a compact manifold $X$ and $g$ a bundle like metric. As above $\hat{H}^n (X,\Rh)$ denotes the Hilbert space completion of $\oH^n (X,\Rh)$ with respect to the scalar product (\ref{eq:13}). Let $\phi^t$ be an $\Fh$-compatible flow. Then the linear operators $\phi^{t*}$ on $\oH^n (X,\Rh)$ induce a strongly continuous operator group on $\hat{H}^n (X,\Rh)$. In particular the infinitesimal generator $\Theta$ exists as a closed densely defined operator. On $\oH^n (X,\Rh)$ it agrees with the infinitesimal generator in the Fr\'echet topology defined earlier. There exists $\omega > 0$ such that the spectrum of $\Theta$ lies in $- \omega \le \RRe s \le \omega$. If the operators $\phi^{t*}$ are orthogonal then $T = - i \Theta$ is a self-adjoint operator on $\hat{H}^n (X,\Rh) \otimes \C$ and we have
\[
\phi^{t*} = \exp t \Theta = \exp it T
\]
in the sense of the functional calculus for (unbounded) self-adjoint operators on Hilbert spaces. 
\end{theorem}

{\bf Sketch of proof} Estimates show that $\| \phi^{t*} \|$ is locally uniformly bounded in $t$ on $\oH^n (X,\Rh)$. Approximating $h \in \hat{H}^n (X,\Rh)$ by $h_{\nu} \in \oH^n (X,\Rh)$ one now shows as in the proof of the Riemann--Lebesgue lemma that the function $t \mapsto \phi^{t*} h$ is continuous at zero, hence everywhere. Thus $\phi^{t*}$ defines a strongly continuous group on $\hat{H}^n (X,\Rh)$. The remaining assertions follow from semigroup theory \cite{DSch}, Ch. VIII, XII, and in particular from the theorem of Stone.

We now combine theorems \ref{t4}, \ref{t7} and \ref{t8} to obtain the following corollary:

\begin{cor}
  \label{t9}
Let $X$ be a compact $3$-manifold with a foliation $\Fh$ by surfaces having a dense leaf. Let $\phi^t$ be a non-degenerate $\Fh$-compatible flow which is everywhere transversal to $\Fh$. Assume that $\phi^t$ is conformal as in (\ref{eq:18}) with respect to a metric $g$ on $T\Fh$. Then $\Theta$ has pure point spectrum $\Sp^1 (\Theta)$ on $\hat{H}^1 (X,\Rh)$ which is discrete in $\R$ and we have the following equalities of distributions on $\R$:
\begin{eqnarray}
  \label{eq:22}
  \lefteqn{\sum^2_{i=0} (-1)^i \Tr (\phi^* \tei \oH^i (X, \Rh)) = 1 - \sum_{\rho \in \Sp^1 (\Theta)} e^{t\rho} + e^{t\alpha}} \\
& = & \chi_{\Co} (\Fh , \mu) \delta_0 + \sum_{\gamma} l (\gamma) \sum_{k \in \Z \ohne 0} \varepsilon_{\gamma} (k) \delta_{kl (\gamma)} \; . \nonumber
\end{eqnarray}
In the sum the $\rho$'s appear with their geometric multiplicities. All $\rho \in \Sp^1 (\Theta)$ have $\RRe \rho = \frac{\alpha}{2}$.
\end{cor}

{\bf Remarks} 1) Here $e^{t\rho} , e^{t \alpha}$ are viewed as distributions so that evaluated on a test function $\varphi \in \Dh (\R)$ the formula reads:
\begin{equation}
  \label{eq:23}
\Phi (0) - \sum\limits_{\rho \in \Sp^1 (\Theta)} \Phi (\rho) + \Phi (\alpha) = \chi_{\Co} (\Fh , \mu) \varphi (0) + \sum\limits_{\gamma} l (\gamma) \sum\limits_{k \in \Z \ohne 0} \varepsilon_{\gamma} (k) \varphi (k l (\gamma)) \; .
\end{equation}
Here we have put
\[
\Phi (s) = \int_{\R} e^{ts} \varphi (t) \, dt \; .
\]
2) Actually the conditions of the corollary force $\alpha = 0$ i.e. the flow must be isometric with respect to $g$. We have chosen to leave the $\alpha$ in the fomulation since there are good reasons to expect the corollary to generalize to more general phase spaces $X$ than manifolds, where $\alpha \neq 0$ becomes possible i.e. to Sullivan's generalized solenoids. More on this in the next section.\\
3) One can show that the group generated by the lengths of closed orbits is a finitely generated subgroup of $\R$ under the assumptions of the corollary. In order to achieve an infinitely generated group the flow must have fixed points.

\begin{proofof}
  {\bf \ref{t9}}
By \ref{t4}, \ref{t7} we only need to show the equation
\begin{equation}
  \label{eq:24}
  \Tr (\phi^{t*} \tei \oH^1 (X,\Rh)) = \sum_{\rho \in \Sp^1 (\Theta)} e^{t \rho}
\end{equation}
and the assertions about the spectrum of $\Theta$. As in the proof of \ref{t4} one sees that on $\hat{H}^1 (X,\Rh)$ we have
\[
( \phi^{t*} h , \phi^{t*} h') = e^{\alpha t} (h , h') \; .
\]
Hence $e^{-\frac{\alpha}{2} t} \phi^{t*}$ is orthogonal and by the theorem of Stone
\[
T = -iS
\]
is selfadjoint on $\hat{H}^1 (X,\Rh) \otimes \C$, if $\Theta = \frac{\alpha}{2} + S$. Moreover
\[
e^{-\frac{\alpha}{2} t} \phi^{t*} = \exp it T \; ,
\]
so that
\begin{equation}
  \label{eq:25}
  \phi^{t*} = \exp t \Theta \; .
\end{equation}
In \cite{DS2} proof of 2.6, for isometric flows the relation
\[
-\Theta^2 = \Delta^1 \, |_{\ker \Delta^1_{\Fh}}
\]
was shown. Using the spectral theory of the ordinary Laplacian $\Delta^1$ on $1$-forms it follows that $\Theta$ has pure point spectrum with finite multiplicities on $\oH^1 (X,\Rh) \cong \ker \Delta^1_{\Fh}$ and that $\Sp^1 (\Theta)$ is discrete in $\R$. Alternatively, without knowing $\alpha = 0$, that proof gives:
\[
- \left( \Theta - \frac{\alpha}{2} \right)^2 = \Delta^1 \, |_{\ker \Delta^1_{\Fh}} \; .
\]
This also implies the assertion on the spectrum of $\Theta$ on $\hat{H}^1 (X,\Rh)$. 
\end{proofof}

\begin{remarksnn}
{\bf a} \quad Writing the eigenvalues $\rho$ of $\Theta$ on $\hat{H}^1 (X, \Rh)$ in the form $\rho = \frac{\alpha}{2} + ir$, the proof shows that the numbers $r^2$ are contained in the spectrum of the ordinary Laplacian $\Delta^1 = \Delta^1_{\tilde{g}}$ on $1$-forms on $X$. Here $\tilde{g}$ is any extension of the metric $g$ on $T\Fh$ to a bundle-like metric on $TX$. One may wonder whether there is a corresponding statement for the zeroes of $\hzeta_K (s)$: Are the squares of their imaginary parts contained in the spectrum of suitable Laplace--Bertrami operators on $1$-forms?\\
{\bf b} \quad In more general situations where $\Theta$ may not have a pure point spectrum on $\hat{H}^1 (X,\Rh)$ but where $e^{-\frac{\alpha}{2}} \phi^{t*}$ is still orthogonal, we obtain:
\[
\langle \Tr (\phi^* \tei \oH^1 (X,\Rh)) , \varphi \rangle = \sum_{\rho \in \Sp^1 (\Theta)_{\mathrm{point}}} \Phi (\rho) + \int^{\frac{\alpha}{2} + i \infty}_{\frac{\alpha}{2} - i \infty} \Phi (\lambda) m (\lambda) \, d \lambda
\]
where $m (\lambda) \ge 0$ is the spectral density function of the continuous part of the spectrum of $\Theta$. 
\end{remarksnn}

\section{Comparison with the ``explicit formulas'' in analytic number theory}
Consider a number field $K / \Q$. 
Formula (\ref{eq:20}) restricted to $\R^{> 0}$ implies the following equality of distributions on $\R^{> 0}$:
\begin{eqnarray}
  \label{eq:26}
  \lefteqn{\sum^2_{i=0} (-1)^i \Tr (\phi^* \tei H^i_{\dyn} (\overline{\spec \eo} , \Rh)) = 1 - \sum_{\hat{\zeta}_K (\rho) = 0} e^{t\rho} + e^t } \\
 & = & \sum_{\ep \nmid \infty} \log N\ep \sum^{\infty}_{k=1} \delta_{k \log N \ep} + \sum_{\ep \tei \infty} (1 - e^{\kappa_{\ep} t})^{-1} \; . \nonumber
\end{eqnarray}
This fits rather nicely with formula (\ref{eq:22}) and suggests the following analogies:
\begin{center}
  \begin{tabular}{llp{8cm}}
$\spec \eo_K \cup \{ \ep \tei \infty \}$ & \quad $\ent$ \quad & $3$-dimensional dynamical system $(X , \phi^t)$ with a one-codimensional foliation $\Fh$ satisfying the conditions of conjecture \ref{t5} \\
finite place $\ep$ & \quad $\ent$ \quad & closed orbit $\gamma = \gamma_{\ep}$ not contained in a leaf and hence transversal to $\Fh$ such that $l (\gamma_{\ep}) = \log N \ep$ and $\varepsilon_{\gamma_{\ep}} (k) = 1$ for all $k \ge 1$.\\
infinite place $\ep$ & \quad $\ent$ \quad & fixed point $x_{\ep}$ such that $\kappa_{x_{\ep}} = \kappa_{\ep}$ and \nolinebreak $\varepsilon_{x_{\ep}} = 1$. \\
$H^i_{\dyn} (\overline{\spec \eo} , \Rh)$ & \quad $\ent$ \quad & $\oH^i (X , \Rh)$
  \end{tabular}
\end{center}

In order to understand number theory more deeply in geometric terms it would be very desirable to find a system $(X , \phi^t , \Fh)$ which actually realizes this correspondence. For this the class of compact $3$-manifolds as phase spaces has to be generalized as will become clear from the following discussion.

Let us compare formula (\ref{eq:4}) on all of $\R$ with formula (\ref{eq:22}) in Corollary \ref{t9}. This corollary is the best result yet on the dynamical side but still only a first step since it does not allow for fixed points which as we have seen must be expected for dynamical systems of relevance for number fields.

Ignoring the contributions $W_{\ep}$ from the infinite places for the moment we are suggested that
\begin{equation}
  \label{eq:27}
  - \log |d_{K/ \Q}| \ent \chi_{Co} (\Fh , \mu) \; .
\end{equation}
There are two nice points about this analogy. Firstly there is the following well known fact due to Connes: 

\begin{fact}
  \label{t10}
Let $\Fh$ be a foliation of a compact $3$-manifold by surfaces such that the union of the compact leaves has $\mu$-measure zero, then
\[
\chi_{Co} (\Fh , \mu) \le 0 \; .
\]
\end{fact}
Namely the non-compact leaves are known to be complete in the induced metric. Hence they carry no non-zero harmonic $L^2$-functions, so that Connes' $0$-th Betti number $\beta_0 (\Fh , \mu) = 0$. Since $\beta_2 (\Fh , \mu) = \beta_0 (\Fh , \mu)$ it follows that
\begin{eqnarray*}
  \chi_{Co} (\Fh , \mu) & = & \beta_0 (\Fh , \mu) - \beta_1 (\Fh , \mu) + \beta_2 (\Fh , \mu) \\
& = & - \beta_1 (\Fh , \mu) \le 0 \; .
\end{eqnarray*}
The reader will have noticed that in accordance with \ref{t10} the left hand side of (\ref{eq:27}) is negative as well:
\[
- \log |d_{K / \Q}| \le 0 \quad \mbox{for all} \;  K / \Q \; .
\]
The second nice point about (\ref{eq:27}) is this. The bundle-like metric $g$ which we have chosen for the definition of $\Delta_{\Fh}$ and of $\chi_{Co} (\Fh , \mu)$ induces a holomorphic structure of $\Fh$ \cite{MS}, Lemma A\,3.1. The space $X$ is therefore foliated by Riemann surfaces. Let $\chi_{Co} (\Fh , \Oh , \mu)$ denote the holomorphic Connes Euler characteristic of $\Fh$ defined using $\Delta_{\overline{\partial}}$-harmonic forms on the leaves instead of $\Delta$-harmonic ones. According to Connes' Riemann--Roch Theorem \cite{MS} Cor. A.\,2.3, Lemma A\,3.3 we have:
\[
\chi_{Co} (\Fh , \Oh , \mu) = \halb \chi_{Co} (\Fh , \mu) \; .
\]
Therefore $\chi_{Co} (\Fh , \Oh , \mu)$ corresponds to $-\log \sqrt{|d_{K / \Q}|}$. 

Now, for completely different reasons this number is defined in Arakelov theory as the Arakelov Euler characteristic of $\overline{\spec \eo_K}$:
\begin{equation}
  \label{eq:28}
  \chi_{Ar} (\Oh_{\overline{\spec \eo_K}}) = - \log {\textstyle \sqrt{|d_{K / \Q}|}} \; .
\end{equation}
See \cite{N} for example. Thus we see that
\[
\chi_{Ar} (\Oh_{\overline{\spec \eo_K}}) \quad \mbox{corresponds to} \; \chi_{Co} (\Fh , \Oh , \mu) \; .
\]
It would be very desirable of course to understand Arakelov Euler characteristics in higher dimensions even conjecturally in terms of Connes' holomorphic Euler characteristics. Note however that Connes' Riemann--Roch theorem in higher dimensions does not involve the $R$-genus appearing in the Arakelov Riemann--Roch theorem. The ideas of Bismut \cite{Bi} may be relevant in this connection. He interpretes the $R$-genus in a natural way via the geometry of loop spaces. 

Further comparison of formulas (\ref{eq:4}) and (\ref{eq:22}) shows that in a dynamical system corresponding to number theory we must have $\alpha = 1$. This means that the flow $\phi^{t*}$ would act by multiplication  with $e^t$ on the one-dimensional space $\oH^2_{\Fh} (X)$. As explained before this would be the case if $\phi^t$ were conformal on $T\Fh$ with factor $e^t$:
\begin{equation}
  \label{eq:29}
  g (T_x \phi^t (v) , T_x \phi^t (w)) = e^t g (v,w) \quad \mbox{for all} \; v,w \in T_x \Fh \; .
\end{equation}
However as mentioned before, this is not possible in the manifold setting of corollary \ref{t9} which actually implies $\alpha = 0$. \\
An equally important difference between formulas (\ref{eq:4}) and (\ref{eq:22}) is between the coefficients of $\delta_{kl (\gamma)}$ and of $\delta_{k \log N\ep}$ for $k \le -1$. In the first case it is $\pm 1$ whereas in the second it is $N\ep^k = e^{k \log N \ep}$ which corresponds to $e^{kl (\gamma)}$. 

Thus it becomes vital to find phase spaces $X$ more general than manifolds for which the analogue of corollary \ref{t9} holds and where $\alpha \neq 0$ and in particular $\alpha = 1$ becomes possible. In the new context the term $\varepsilon_{\gamma} (k) \delta_{k l (\gamma)}$ for $k \le -1$ in formula (\ref{eq:22}) should become $\varepsilon_{\gamma} (k) e^{\alpha k l (\gamma)} \delta_{kl (\gamma)}$. The next section is devoted to a discussion of certain laminated spaces which we propose as possible candidates for this goal.

\section{Remarks on dynamical Lefschetz trace formulas on laminated spaces}
In this section we extend the previous discussion to more general phase spaces than manifolds. The class of spaces we have in mind are the smooth solenoids \cite{Su} or in other words foliated spaces with totally disconnected transversals in the sense of \cite{MS}. 

\begin{punkt} \label{t11} \rm
We start by recalling the definition of a foliated space with $d$-dimensional leaves in the sense of \cite{MS} Cap. II. Consider a separable metrizable topological space $X$ with a covering by open sets $U_i$ and homeomorphisms $\varphi_i : U_i \to F_i \times T_i$ with $F_i$ open in $\R^d$. This atlas should have the following property\\
{\bf a} \quad The transition functions between different charts $\varphi_i$ and $\varphi_j$ have the following form locally where $x$ is the Euklidean and $y$ the transversal coordinate
\[
\varphi_j \verk \varphi^{-1}_i (x,y) = (f_{ij} (x,y) , g_{ij} (y)) \; .
\]
{\bf b} \quad All partial derivatives $D^{\alpha}_x f_{ij} (x,y)$ exist and are continuous as functions of $x$ and $y$.

Two atlases of this kind are equivalent if their union is an atlas as well. A foliated space with $d$-dimensional leaves is a space $X$ as above together with an equivalence class of atlases. The coordinate change locally transforms the level surface or plaque $\varphi^{-1}_i (F_i \times \{ y \} )$ to the level surface $\varphi^{-1}_j (F_j \times \{ y' \})$ where $y' = g_{ij} (y)$. Hence the level surfaces glue to maximal connected sets with a smooth $d$-dimensional manifold structure. These are the leaves $F$ of the foliated space $X$. Usually we indicate the foliated structure on a topological space $X$ by another symbol, e.g. $\Fh$ which we also use to denote the set of leaves.

A (generalized) solenoid of dimension $a$ is a separable metrizable space $X$ which is locally homeomorphic to a space of the form $L_i \times T_i$ where $L_i$ is an open subset of $\R^a$ and $T_i$ is totally disconnected. For such spaces the transition functions between charts automatically satisfy condition {\bf a} since continuous functions from connected subsets of $\R^a$ into a totally disconnected space are constant. If $X$ carries a smooth structure, i.e. if we are given an atlas of charts whose coordinate changes satisfy condition {\bf b}, then $X$ becomes a foliated space with $a$-dimensional leaves. For this foliation $\Lh$ the leaves $L$ of $\Lh$ are the path-connected components of $X$.

The two most prominent places in mathematics where (smooth) solenoids occur naturally are in number theory e.g. as adelic points of algebraic groups and in the theory of dynamical systems as attractors.

The classical solenoid
\[
\Sa^1_p = \R \times_{\Z} \Z_p = \varprojlim ( \ldots \xrightarrow{p} \R / \Z \xrightarrow{p} \R / \Z \xrightarrow{p} \ldots)
\]
is an example of a compact connected one-dimensional smooth solenoid with dense leaves diffeomorphic to the real line.

For a foliated space $(X, \Fh)$ let $T\Fh$ be the tangent bundle in the sense of \cite{MS} p. 43. For a point $x \in X$, the fibre $T_x \Fh$ is the usual tangent space of the leaf through the point $x$. Morphisms between foliated spaces are continuous maps which induce smooth maps between the leaves. They induce morphisms of tangent bundles.

For a smooth solenoid $(X, \Lh)$ we set $TX = T\Lh$. By definition, a Riemannian metric on $X$ is one on $TX$.
\end{punkt}

\begin{punkt} \label{t12}
\rm We now turn to cohomology. For a foliated space $(X,\Fh)$ of leaf dimension $d$ consider the sheaf $\Rh = \Rh_{\Fh}$ of real valued continuous functions which are locally constant on the leaves. There is a natural de~Rham resolution 
\[
0 \longrightarrow \Rh \longrightarrow \Ah^0_{\Fh} \xrightarrow{d_{\Fh}} \Ah^1_{\Fh} \xrightarrow{d_{\Fh}} \ldots
\]
by the fine sheaves of differential forms along the leaves. Explicitly
\[
\Ah^i_{\Fh} (U) = \Gamma (U , \Lambda^i T^* \Fh)
\]
for every open subset $U$ of $X$. Here sections are by definition continuous and smooth on the leaves. Hence we have
\[
H^i (X, \Rh) = \frac{\Ker (d_{\Fh} : \Ah^i_{\Fh} (X) \to \Ah^{i+1}_{\Fh} (X))}{\Imm (d_{\Fh} : \Ah^{i-1}_{\Fh} (X) \to \Ah^i_{\Fh} (X))} \; .
\]
As before, one also considers the maximal Hausdorff quotient $\oH^i (X,\Rh)$ of this cohomology, obtained by dividing by the closure of $\Imm d_{\Fh}$ in the natural Fr\'echet topology.
\bigskip

\noindent {\bf Warning} A manifold with a (smooth) foliation is also a foliated space. However the sheaves $\Rh$ and $\Ah^i_{\Fh}$ are different in the two contexts: In the first one demands smoothness also in the transversal direction whereas in the second one only wants continuity.
\bigskip

Now consider an $a$-dimensional smooth solenoid $(X, \Lh)$ with a one-codi\-mensional foliation $\Fh$. By this we mean that $X$ is equipped with a further structure $\Fh$ of a foliated space, this time with leaves of dimension $d = a-1$ such that the leaves $L$ of $\Lh$ are foliated by the leaves $F$ of $\Fh$ with $F \subset L$. We denote by $\Fh_L$ the induced one-codimensional foliation of $L$.
\end{punkt}

\begin{punkt}
  \label{t14} \rm
A flow $\phi$ on $X$ is a continuous $\R$-action such that the induced $\R$-actions on the leaves of $\Lh$ are smooth. It is compatible with $\Fh$ if every $\phi^t$ maps leaves of $\Fh$ into leaves of $\Fh$. Thus $(X , \Fh , \phi^t)$ is partitioned into the foliated dynamical systems $(L, \Fh_L , \phi^t \, |_L)$ for $L \in \Lh$. Any $\Fh$-compatible flow $\phi^t$ induces pullback actions $\phi^{t*}$ on $H^{\hullet} (X, \Rh_{\Fh})$ and $\oH^{\hullet} (X, \Rh_{\Fh})$.
\end{punkt}

\begin{punkt}
  \label{t15} \rm
We now state as a working hypotheses a generalization of the conjectured dynamical trace formula \ref{t5}. We allow the phase space to be a smooth solenoid. Moreover we extend the formula to an equality of distributions on $\R^*$ instead of $\R^{>0}$. After checking various compatibilities we state a case where our working hypotheses can be proved and give a number theoretical example.
\end{punkt}

\begin{punkt} {\bf Working hypotheses:}
  \label{t16} \rm
Let $X$ be a compact smooth solenoid with a one-codimensional foliation $\Fh$ and an $\Fh$-compatible flow $\phi$. Assume that the fixed points and the periodic orbits of the flow are non-degenerate. Then there exists a natural definition of a $\Dh' (\R^*)$-valued trace of $\phi^{t*}$ on $\oH^{\hullet} (X, \Rh)$ where $\Rh = \Rh_{\Fh}$ such that in $\Dh' (\R^*)$ we have:
\begin{eqnarray}
  \label{eq:33}
\hspace*{0.5cm}  \lefteqn{\sum^{\dim \Fh}_{n=0} (-1)^n \Tr (\phi^* \tei \oH^n (X , \Rh)) = }\\
& & \sum_{\gamma} l (\gamma) \left( \sum_{k \ge 1} \varepsilon_{\gamma} (k) \delta_{k l (\gamma)} + \sum_{k \le -1} \varepsilon_{\gamma} (|k|) \det (-T_x \phi^{kl (\gamma)} \tei T_x \Fh) \delta_{kl (\gamma)} \right) \nonumber \\
& & + \sum_x W_x \; . \nonumber
\end{eqnarray}
Here $\gamma$ runs over the closed orbits not contained in a leaf and in the sums over $k$'s any point $x \in \gamma$ can be chosen. The second sum runs over the fixed points $x$ of the flow. The distributions $W_x$ on $\R^*$ are given by:
\[
W_x \, |_{\R^{> 0}} = \varepsilon_x \, |1 - e^{\kappa_x t}|^{-1}
\]
and
\[
W_x \, |_{\R^{< 0}} = \varepsilon_x \det (-T_x \phi^t \tei T_x \Fh) \, |1 - e^{\kappa_x |t|}|^{-1} \; .
\]
\end{punkt}

\begin{remarks} \label{t17}
  \rm {\bf 0)} It may actually be better to use a version of foliation cohomology where transversally forms are only supposed to be locally $L^2$ instead of being continuous. With such a foliation $L^2$-cohomology, Leichtnam \cite{Lei} has proved certain fixed point free cases of the working hypothesis by generalizing the techniques of \cite{AK1}, \cite{AK2} to the solenoidal setting.\\
{\bf 1)} In the situation described in \ref{t18} below the working hypotheses can be proved, c.f. Theorem \ref{t19}. In those cases there are no fixed points, only closed orbits. Thus Theorem \ref{t19} dictated only the coefficients of $\delta_{kl (\gamma)}$ for $k \in \Z \ohne 0$, but not the contributions $W_x$ from the fixed points. \\
{\bf 2)} The coefficients of $\delta_{kl (\gamma)}$ for $k \in \Z \ohne 0$ can be written in a uniform way as follows. They are equal to:
\begin{equation}
  \label{eq:34}
  \frac{\det (1 - T_x \phi^{kl (\gamma)} \tei T_x \Fh)}{|\det (1 - T_x \phi^{|k| l (\gamma)} \tei T_x X / \R  Y_{\phi , x})|} = \frac{\det (1 - T_x \phi^{kl (\gamma)} \tei T_x \Fh)}{| \det (1 - T_x \phi^{|k| l (\gamma)} \tei T_x \Fh)|} \; .
\end{equation}
Here $x$ is any point on $\gamma$. Namely, for $k \ge 1$ this equals $\varepsilon_{\gamma} (k)$ whereas for $k \le -1$ we obtain
\begin{equation}
  \label{eq:35}
  \varepsilon_{\gamma} (|k|) \det (- T_x \phi^{kl (\gamma)} \tei T_x \Fh) = \varepsilon_{\gamma} (k) \, |\det (T_x \phi^{kl (\gamma)} \tei T_x \Fh)| \; .
\end{equation}
The expression on the left hand side of (\ref{eq:34}) motivated our conjecture about the contributions on $\R^*$ from the fixed points $x$. Since $Y_{\phi , x} = 0$, they should be given by:
\[
\frac{\det (1 - T_x \phi^t \tei T_x \Fh)}{| \det (1 - T_x \phi^{|t|} \tei T_x X)|} \overset{!}{=} W_x \; .
\]
{\bf 3)} One can prove that in the manifold setting of theorem \ref{t7} we have
\[
|\det (T_x \phi^{kl (\gamma)} \tei T_x \Fh)| = 1 \; .
\]
By (\ref{eq:35}), our working hypotheses \ref{t16} is therefore compatible with formula (\ref{eq:21}). Compatibility with conjecture \ref{t5} is clear.\\
{\bf 4)} We will see below that in our new context metrics $g$ on $T\Fh$ can exist for which the flow has the conformal behaviour (\ref{eq:29}). Assuming we are in such a situation and that $\Fh$ is $2$-dimensional, we have:
\[
|\det (T_x \phi^{kl (\gamma)} \tei T_x \Fh)| = e^{kl (\gamma)} \quad \mbox{for} \; x \in \gamma , k \in \Z
\]
and
\[
|\det (T_x \phi^t \tei T_x \Fh)| = e^t \quad \mbox{for a fixed point} \; x \; .
\]
In the latter case, we even have by continuity:
\[
\det (T_x \phi^t \tei T_x \Fh) = e^t \; ,
\]
the determinant being positive for $t = 0$. Hence by (\ref{eq:35}) the conjectured formula (\ref{eq:33}) reads as follows in this case:
\begin{eqnarray}
  \label{eq:36}
  \lefteqn{\sum^2_{n=0} (-1)^n \Tr (\phi^* \tei \oH^n (X, \Rh))} \\
& = & \sum_{\gamma} l (\gamma) \left( \sum_{k \ge 1} \varepsilon_{\gamma} (k) \delta_{kl (\gamma)} + \sum_{k \le -1} \varepsilon_{\gamma} (k) e^{kl (\gamma)} \delta_{kl (\gamma)} \right) + \sum_x W_x \; . \nonumber
\end{eqnarray}
Here:
\[
W_x \, |_{\R^{> 0}} = \varepsilon_x \, |1 - e^{\kappa_x t}|^{-1}
\]
and
\[
W_x \, |_{\R^{< 0}} = \varepsilon_x e^t \, |1 - e^{\kappa_x |t|}|^{-1} \; .
\]
This fits perfectly with the explicit formula (\ref{eq:5}) if all $\varepsilon_{\gamma_{\ep}} (k) = 1$ and $\varepsilon_{x_{\ep}} = 1$. Namely if $l (\gamma_{\ep}) = \log N \ep$ for $\ep \nmid \infty$ and $\kappa_{x_{\ep}} = \kappa_{\ep}$ for $\ep \tei \infty$, then we have:
\[
e^{kl (\gamma_{\ep})} = e^{k \log N\ep} = N \ep^k \quad \mbox{for finite places} \; \ep
\]
and
\[
W_{x_{\ep}} = W_{\ep} \quad \mbox{on} \; \R^* \; \mbox{for the infinite places} \; \ep \; .
\]
{\bf 5)} In the setting of the preceeding remark the automorphisms
\[
e^{-\frac{k}{2} l (\gamma)} T_x \phi^{kl (\gamma)} \quad \mbox{of} \; T_x \Fh \quad \mbox{for} \; x \in \gamma
\]
respectively
\[
e^{-\frac{t}{2}} T_x \phi^t \quad \mbox{of} \; T_x \Fh \quad \mbox{for a fixed point} \; x
\]
are orthogonal automorphisms. For a real $2 \times 2$ orthogonal determinant $O$ with $\det O = -1$ we have:
\[
\det (1 - uO) = 1 - u^2 \; .
\]
The condition $\varepsilon_{\gamma} (k) = +1$ therefore implies that $\det (T_x \phi^{kl (\gamma)} \tei T_x \Fh)$ is positive for $k \ge 1$ and hence for all $k \in \Z$. The converse is also true. For a fixed point we have already seen directly that $\det (T_x \phi^t \tei T_x \Fh)$ is positive for all $t \in \R$. Hence we have the following information. 

{\bf Fact} In the situation of the preceeding remark, $\varepsilon_k (\gamma) = +1$ for all $k \in \Z \ohne 0$ if and only if on $T_x \Fh$ we have:
\[
T_x \phi^{kl (\gamma)} = e^{\frac{k}{2} l (\gamma)} \cdot O_k \quad \mbox{for} \; O_k \in \SO (T_x \Fh) \; .
\]
For fixed points, $\varepsilon_x = 1$ is automatic and we have:
\[
T_x \phi^t = e^{\frac{t}{2}} O_t \quad \mbox{for} \; O_t \in \SO (T_x \Fh) \; .
\]

In the number theoretical case the eigenvalues of $T_x \phi^{\log N\ep}$ on $T_x \Fh$ for $x \in \gamma_{\ep}$ would therefore be complex conjugate numbers of absolute value $N\ep^{1/2}$. If they are real then $T_x \phi^{\log N\ep}$ would simply be mutliplication by $\pm N\ep^{1/2}$. If not, the situation would be more interesting. Are the eigenvalues Weil numbers (of weight $1$)? If yes there would be some elliptic curve over $\eo_K / \ep$ involved by Tate--Honda theory. \\
{\bf 6)} It would of course be very desirable to extend the hypotheses \ref{t16} to a conjectured equality of distributions on all of $\R$. By theorem \ref{t7} we expect one contribution of the form
\[
\chi_{\Co} (\Fh , \mu) \cdot \delta_0 \; .
\]
The analogy with number theory suggests that there will also be somewhat complicated contributions from the fixed points in terms of principal values which are hard to guess at the moment. After all, even the simpler conjecture \ref{t5} has not yet been verified in the presence of fixed points!\\
{\bf 7)} If there does exist a foliated dynamical system attached to $\overline{\spec \eo_K}$ with the properties dictated by our considerations we would expect in particular that for a preferred transverse measure $\mu$ we have:
\[
\chi_{\Co} (\Fh , \mu) = - \log |d_{K / \Q}| \; .
\]
This gives some information on the space $X$ with its $\Fh$-foliation. If $K / \Q$ is ramified at some finite place i.e. if $d_{K / \Q} \neq \pm 1$ then $\chi_{\Co} (\Fh , \mu) < 0$. Now, since $\oH^2 (X, \Fh)$ must be one-dimensional, it follows that
\[
\chi_{\Co} (\Fh , \nu) < 0 \quad \mbox{for all non-trivial transverse measures} \; \nu \; .
\]
Hence by a result of Candel \cite{Ca} there is a Riemannian metric on $T \Fh$, such that every $\Fh$-leaf has constant curvature $-1$. Moreover $(X , \Fh)$ is isomorphic to
\[
\Oh (H , X) / \PSO (2) \; .
\]
Here $\Oh (H , X)$ is the space of conformal covering maps $u : H \to N$ as $N$ runs through the leaves of $\Fh$ with the compact open topology. See \cite{Ca} for details.

In the unramified case, $|d_{K / \Q}| = 1$ we must have $\chi_{\Co} (\Fh , \nu) = 0$ for all transverse measures by the above argument. Hence there is an $\Fh$-leaf which is either a plane, a torus or a cylinder c.f. \cite{Ca}.
\end{remarks}

\begin{punkt}
  \label{t18}
\rm In this final section we describe a simple case where the working hypothesis \ref{t16} can be proved.

Consider an unramified covering $f : M \to M$ of a compact connected orientable $d$-dimensional manifold $M$. We set
\[
\bar{M} = \lim_{\leftarrow} ( \ldots \xrightarrow{f} M \xrightarrow{f} M \to \ldots) \; .
\]
Then $\bar{M}$ is a compact topological space equipped with the shift automorphism $\bar{f}$ induced by $f$. It can be given the structure of a smooth solenoid as follows. Let $\tilde{M}$ be the universal covering of $M$. For $i \in \Z$ there exists a Galois covering
\[
p_i : \tilde{M} \longrightarrow M
\]
with Galois group $\Gamma_i$ such that $p_i = p_{i+1} \verk f$ for all $i$. Hence we have inclusions:
\[
\ldots \subset \Gamma_{i+1} \subset \Gamma_i \subset \ldots \subset \Gamma_0 =: \Gamma \cong \pi_1 (M , x_0) \; .
\]
Writing the operation of $\Gamma$ on $\tilde{X}$ from the right, we get commutative diagrams for $i \ge 0$:
\[
\begin{CD}
  \tilde{M} \times_{\Gamma} (\Gamma / \Gamma_{i+1}) @= \tilde{M} / \Gamma_{i+1} @>{\overset{p_{i+1}}{\sim}}>> M \\
@VV{\id \times \proj}V  @VV{\proj}V @VV{f}V \\
\tilde{M} \times_{\Gamma} (\Gamma / \Gamma_i) @= \tilde{M} / \Gamma_i @>{\overset{p_i}{\sim}}>> M
\end{CD}
\]
It follows that
\begin{equation}
  \label{eq:37}
  \tilde{M} \times_{\Gamma} \bar{\Gamma} \silo \bar{M}
\end{equation}
where $\bar{\Gamma}$ is the pro-finite set with $\Gamma$-operation:
\[
\bar{\Gamma} = \lim_{\leftarrow} \Gamma / \Gamma_i \; .
\]
The isomorphism (\ref{eq:37}) induces on $\bar{M}$ the structure of a smooth solenoid with respect to which $\bar{f}§$ becomes leafwise smooth.

Fix a positive number $l > 0$ and let $\Lambda = l \Z \subset \R$ act on $\bar{M}$ as follows: $\lambda = l \nu$ acts by $\bar{f}^{\nu}$. Define a right action of $\Lambda$ on $\bar{M} \times \R$ by the formula
\[
(m,t) \cdot \lambda = (- \lambda \cdot m , t + \lambda) = (\bar{f}^{-\lambda / l} (m) , t + \lambda) \; .
\]
The suspension:
\[
X = \bar{M} \times_{\Lambda} \R = (\tilde{M} \times_{\Gamma} \overline{\Gamma}) \times_{\Lambda} \R
\]
is an $a = d+1$-dimensional smooth solenoid with a one-codimensional foliation $\Fh$ as in \ref{t12}. The leaves of $\Fh$ are the images in $X$ of the manifolds $\tilde{M} \times \{ \overline{\gamma} \} \times \{ t \}$ for $\overline{\gamma} \in \overline{\Gamma}$ and $t \in \R$. Translation in the $\R$-variable
\[
\phi^t [m,t'] = [m,t+t']
\]
defines an $\Fh$-compatible flow $\phi$ on $X$ which is everywhere transverse to the leaves of $\Fh$ and in particular has no fixed points.

The map
\[
\gamma \longmapsto \gamma_M = \gamma \cap (\bar{M} \times_{\Lambda} \Lambda)
\]
gives a bijection between the closed orbits $\gamma$ of the flow on $X$ and the finite orbits $\gamma_M$ of the $\bar{f}$- or $\Lambda$-action. These in turn are in bijection with the finite orbits of the original $f$-action on $M$. We have:
\[
l (\gamma) = |\gamma_M| l \; .
\]
\end{punkt}

\begin{theorem}
  \label{t19}
In the situation of \ref{t18} assume that all periodic orbits of $\phi$ are non-degenerate. Let $\Sp^n (\Theta)$ denote the set of eigenvalues with their algebraic multiplicities of the infinitesimal generator $\Theta$ of $\phi^{t*}$ on $\oH^n (X,\Rh)$. Then the trace
\[
\Tr (\phi^* \tei \oH^n (X,\Rh)) := \sum_{\lambda \in \Sp^n (\Theta)} e^{t\Theta}
\]
defines a distribution on $\R$ and the following formula holds true in $\Dh' (\R)$:
\begin{eqnarray*}
  \lefteqn{\sum^{\dim \Fh}_{n=0} (-1)^n \Tr (\phi^* \tei \oH^n (X,\Rh)) = \chi_{\Co} (\Fh , \mu) \cdot \delta_0 \; + }\\
&& \sum_{\gamma} l (\gamma) \left( \sum_{k \ge 1} \varepsilon_{\gamma} (k) \delta_{kl (\gamma)} + \sum_{k \le -1} \varepsilon_{\gamma} (|k|) \det (-T_x \phi^{kl (\gamma)} \tei T_x \Fh) \delta_{kl (\gamma)} \right) \; .
\end{eqnarray*}
Here $\gamma$ runs over the closed orbits of $\phi$ and in the sum over $k$'s any point $x \in \gamma$ can be chosen. Moreover $\chi_{\Co} (\Fh , \mu)$ is the Connes' Euler characteristic of $\Fh$ with respect to a certain canonical transverse measure $\mu$. Finally we have the formula:
\[
\chi_{\Co} (\Fh , \mu) = \chi (M) \cdot l \; .
\]
\end{theorem}

A more general result is proved in \cite{Lei}.

\begin{example}
  Let $E / \F_p$ be an ordinary elliptic curve over $\F_p$ and let $\C / \Gamma$ be a lift of $E$ to a complex elliptic curve with $CM$ by the ring of integers $\eo_K$ in an imaginary quadratic field $K$. Assume that the Frobenius endomorphism of $E$ corresponds to the prime element $\pi$ in $\eo_K$. Then $\pi$ is split, $\pi \bar{\pi} = p$ and for any embedding $\Q_l \subset \C , l \neq p$ the pairs
\[
(H^*_{\et} (E \otimes \bar{\F}_p , \Q_l) \otimes \C , \Frob^*) \quad \mbox{and} \quad (H^* (\C / \Gamma , \C) , \pi^*)
\]
are isomorphic. Setting $M = \C / \Gamma , f = \pi$ we are in the situation of \ref{t18} and we find:
\[
X = (\C \times_{\Gamma} T_{\pi} \Gamma) \times_{\Lambda} \R \; .
\]
Here
\[
T_{\pi} \Gamma = \lim_{\leftarrow} \Gamma / \pi^i \Gamma \cong \Z_p
\]
is the $\pi$-adic Tate module of $\C / \Gamma$. It is isomorphic to the $p$-adic Tate module of $E$.
\end{example}

Setting $l = \log p$, so that $\Lambda = (\log p) \Z$ and passing to multiplicative time, $X$ becomes isomorphic to
\[
X \cong (\C \times_{\Gamma} T_{\pi} \Gamma) \times_{p^{\Z}} \R^*_+
\]
which may be a more natural way to write $X$. Note that $p^{\nu}$ acts on $\C \times_{\Gamma} T_{\pi} \Gamma$ by diagonal multiplication with $\pi^{\nu}$. It turns out that the right hand side of the dynamical Lefschetz trace formula established in theorem \ref{t19} equals the right hand side in the explicit formulas for $\zeta_E (s)$. Moreover the metric $g$ on $T \Fh$ given by
\[
g_{[z,y,t]} (\xi, \eta) = e^t \RRe (\xi \bar{\eta}) \quad \mbox{for} \; [z,y,t]
\; \mbox{in}\; (\C \times_{\Gamma} T_{\pi} \Gamma) \times_{\Lambda} \R 
\]
satisfies the conformality condition (\ref{eq:18}) for $\alpha = 1$. The construction of $(X , \phi^t)$ that we made for ordinary elliptic curves is misleading however, since it almost never happens that a variety in characteristic $p$ can be lifted to characteristic zero {\it together with its Frobenius endomorphism}. 

Our present dream for the general situation is this: To an algebraic scheme $\eX / \Z$ one should first attach an infinite dimensional dissipative dynamical system, possibly using $\GL_{\infty}$ in some way. The desired dynamical system should then be obtained by passing to the finite dimensional compact global attractor, c.f. \cite{La} Part I. 


\section{Arithmetic topology and dynamical systems}
In the sixties Mazur and Manin pointed out intriguing analogies between prime ideals in number rings and knots in $3$-manifolds. Let us recall some of the relevant ideas. 

In many respects the spectrum of a finite field $\F_q$ behaves like a topological circle. For example its \'etale cohomology with $\Z_l$-coefficients for $l \nmid q$ is isomorphic to $\Z_l$ in degrees $0$ and $1$ and it vanishes in higher degrees. 

Artin and Verdier \cite{AV} have defined an \'etale topology on $\overline{\spec \Z} = \\
\spec \Z \cup \{ \infty \}$. In this topology $\overline{\spec \Z}$ has cohomological dimension three, up to $2$-torsion. 

The product formula
\[
\prod_{p \le \infty} |a|_p = 1 \quad \mbox{for} \; a \in \Q^* \; ,
\]
allows one to view $\overline{\spec \Z}$ as a {\it compact} space: Namely in the function field case the analogue of the product formula is equivalent to the formula
\[
\sum_x \ord_x f = 0
\]
on a {\it proper} curve $X_0 / \F_q$. 
Here $x$ runs over the closed points of $X_0$ and $f \in \F_q (X_0)^*$. 

By a theorem of Minkowski there are no non-trivial extensions of $\Q$ unramified at all places $p \le \infty$. Thus 
\[
\hat{\pi}_1 (\overline{\spec \Z}) = 0 \; .
\]
Hence, by analogy with the Poincar\'e conjecture, Mazur suggested to think of $\overline{\spec \Z}$ as an arithmetic analogue of the $3$-sphere $S^3$. Under this analogy the inclusion:
\[
\spec \F_p \hookrightarrow \overline{\spec \Z}
\]
corresponds to an embedded circle i.e. to a knot.

The analogue of the Alexander polynomial of a knot turns out to be the Iwasawa zeta function. One can make this precise using $p$-adic \'etale cohomology for schemes over $\spec \Z_{(p)}$. 

More generally, for a number field $K$ consider
\[
\overline{\spec \eo_K} = \spec \eo_K \cup \{ \ep \tei \infty \}
\]
together with its Artin--Verdier \'etale topology. Via the inclusion
\[
\spec \eo_K / \ep \hookrightarrow \spec \eo_K
\]
we may imagine a prime ideal $\ep$ as being analogous to a knot in a compact $3$-manifold.

This nice analogy between number theory and three dimensional topology was further extended by Reznikov and Kapranov and baptized Arithmetic Topology. The reader may find dictionaries between the two fields in \cite{Re} and \cite{S}. Further analogies were contributed in \cite{Mo} and \cite{Ra} for example. In particular Ramachandran had the idea that the infinite primes of a number field should correspond to the ends of a non-compact manifold -- the analogue of $\spec \eo_K$. 

The arguments in the preceeding sections also lead to the idea that prime ideals $\ep$ in $\eo_K$ are knots in a $3$-space.
Namely, as we have seen a phase space $(\overline{\spec \eo_K})_{\dyn}$ corresponding to $\overline{\spec \eo_K}$ would be three dimensional. After forgetting the parametrization of the periodic orbits the prime ideals can thus be viewed as knots in a $3$-space.

It appears that the (sheaf) cohomology of $(\overline{\spec \eo_K})_{\dyn}$ with constant coefficients should play the role of an arithmetic as opposed to geometric cohomology theory for $\overline{\spec \eo_K}$. In particular it would have cohomological dimension three. We wish to compare this as yet speculative theory with the $l$-adic cohomology of $\overline{\spec \eo_K}$. More precisely we compare Lefschetz numbers of certain endomorphisms on these cohomologies.

In order to do this we first calculate the Lefschetz number  of an automorphism $\sigma$ of $K$ on the Artin--Verdier \'etale cohomology of $\overline{\spec \eo_K}$.

Next we prove a generalization of Hopf's formula to a formula for the Lefschetz number of an endomorphism of a dynamical system on a manifold.

Assuming our formula applies to $((\overline{\spec \eo_K})_{\dyn} , \phi^t)$ we then obtain an expression for the Lefschetz number of the automorphism on $H^{\hullet}_{\dyn} (\overline{\spec \eo_K} , \R)$ induced by $\sigma$.

As it turns out, the two kinds of Lefschetz numbers agree in all cases, even when generalized to constructible sheaf coefficients. This result  was prompted by a question of B. Mazur on the significance of \'etale Euler characteristics  \cite{D1} in our dynamical picture.


We first recall a version of \'etale cohomology with compact supports of arithmetic schemes that takes into account the fibres at infinity. Using this theory we reformulate the main result of \cite{D1} on $l$-adic Lefschetz numbers of $l$-adic sheaves on arithmetic schemes as a vanishing statement. This formulation was suggested by Faltings \cite{F} in his review of \cite{D1}. We also calculate $l$-adic Lefschetz numbers on arithmetically compactified schemes.

Let the scheme $\Uh / \Z$ be algebraic i.e. separated and of finite type and set $\Uh_{\infty} = \Uh^{\ann}_{\C} / G_{\R}$ where $\Uh_{\C} = \Uh \otimes_{\Z} \C$ and the Galois group $G_{\R}$ of $\R$ acts on $\Uh^{\ann}_{\C}$ by complex conjugation. We give $\Uh_{\infty}$ the quotient topology of $\Uh^{\ann}_{\C}$.

Artin and Verdier \cite{AV} define the \'etale topology on $\oUh = \Uh \amalg \Uh_{\infty}$ as follows. The category of ``open sets'' has objects the pairs \\
$(f : \Uh' \to \Uh, D')$ where $f$ is an \'etale morphism and $D' \subset \Uh'_{\infty}$ is open. The map $f_{\infty} : D' \to \Uh_{\infty}$ induced by $f$ is supposed to be ``unramified'' in the sense that $f_{\infty} (D') \in \Uh (\R)$ if and only if $D' \in \Uh' (\R)$. Note that $\Uh (\R)$ is a closed subset of $\Uh_{\infty}$.

A morphism
\[
(f : \Uh' \to \Uh , D') \longrightarrow (g : \Uh'' \to \Uh , D'')
\]
is a map $\Uh' \to \Uh''$ commuting with the structure maps and such that the induced map $\Uh'_{\infty} \to \Uh''_{\infty}$ carries $D'$ into $D''$. Coverings are the obvious ones.

Pullback defines morphisms of sites:
\[
\Uh_{\et} \overset{j}{\longrightarrow} \oUh_{\et} \overset{i}{\longleftarrow} \Uh_{\infty} \; .
\]
Let $\sim$ denote the corresponding categories of abelian sheaves. One proves that $\oUh^{\sim}_{\et}$ is the mapping cone of the left exact functor
\[
i^* j_* : \tilde{\Uh}_{\et} \longrightarrow \tilde{\Uh}_{\infty} \; .
\]
In particular we have maps $i^!$ and $j_!$ at our disposal. Let us describe the functor $i^* j_*$ explicitely. Let $\alpha : \Uh^{\ann}_{\C} \to \Uh_{\et}$ be the canonical map of sites. Note that $\alpha^* F$ is a $G_{\R}$-sheaf on $\Uh^{\ann}_{\C}$ for every sheaf $F$ on $\Uh_{\et}$. If $\pi : \Uh^{\ann}_{\C} \to \Uh_{\infty}$ denotes the natural projection, define the left exact functor
\[
\pi^{G_{\R}}_* : (\mbox{abelian}\;G_{\R}\mbox{-sheaves on} \; \Uh^{\ann}_{\C}) \longrightarrow \tilde{\Uh}_{\infty}
\]
by
\[
\pi^{G_{\R}}_* (G) (V) = G (\pi^{-1} (V))^{G_{\R}} \; .
\]
Then one can check that
\[
i^* j_* = \pi^{G_{\R}}_* \verk \alpha^* \; .
\]
In particular we see that
\[
i^* R^n j_* = R^n \pi^{G_{\R}}_* \verk \alpha^* \; .
\]
It follows that for $n \ge 1$ the sheaf $R^n j_* F = i_* i^* R^n j_* F$ is $2$-torsion with support on $\Uh (\R) \subset \Uh_{\infty}$ as stated in \cite{AV}. In particular
\begin{equation}
  \label{eq:38}
  H^n_{\et} (\oUh , j_* F) \longrightarrow H^n_{\et} (\oUh , R j_* F) = H^n_{\et} (\Uh , F)
\end{equation}
is an isomorphism up to $2$-torsion for $n \ge 1$ (and an isomorphism for $n = 0$).

Let us now define cohomology with compact supports for schemes $\eX$ which are proper over $\spec \Z$. We define:
\[
H^n_c (\eX , F) := H^n (\oeX , j_! F) \; .
\]
The distinguished triangle
\[
j_! j^* G^{\hullet} \longrightarrow G^{\hullet} \longrightarrow i_* i^* G^{\hullet} \longrightarrow \ldots
\]
for complexes of sheaves $G^{\hullet}$ on $\oeX$ applied to $G^{\hullet} = Rj_* F$ gives the triangle
\[
j_! F \longrightarrow Rj_* F \longrightarrow i_* i^* Rj_* F \longrightarrow \ldots
\]
From this one gets an exact sequence
\[
\longrightarrow H^n_c (\eX , F) \longrightarrow H^n (\eX , F) \longrightarrow H^n (\eX_{\infty} , i^* Rj_* F) \longrightarrow \ldots
\]
We have $i^* Rj_* = R\pi^{G_{\R}}_* \verk \alpha^*$ and
\begin{eqnarray}
H^n (\eX_{\infty} , i^* Rj_* F) & = & 
  H^n (\eX_{\infty} , R \pi^{G_{\R}}_* (\alpha^* F)) \nonumber \\
& = & H^n (G_{\R} , R\Gamma (\eX^{\ann}_{\C} , \alpha^* F)) \nonumber \\
& = & H^n (\eX^{\ann}_{\R} , \alpha^* F) \; . \label{eq:39}
\end{eqnarray}
Here, for any $G_{\R}$-sheaf $G$ on a complex analytic space $Y$ with a real structure we set:
\[
H^n (Y_{\R} , G) = H^n (G_{\R} , R \Gamma (Y , G)) \; .
\]
In conclusion we get the long exact sequence:
\[
  \longrightarrow H^n_c (\eX , F) \longrightarrow H^n (\eX, F) \longrightarrow H^n (\eX^{\ann}_{\R} , \alpha^* F) \longrightarrow \ldots
\]
An endomorphism $(\sigma , e)$ of the pair $(\eX , F)$ is a morphism $\sigma : \eX \to \eX$ together with a homomorphism $e : \sigma^* F \to F$ of $\Q_l$-sheaves. It induces pullback endomorphisms on $H^n_c (\eX , F)$ and $H^n (\eX , F)$ and on $H^n_c (\eX^{\ann}_{\R} , \alpha^* F)$. For an endomorphism $\varphi$ of a finite dimensional graded vector space $H^{\hullet}$ we will abbreviate the Lefschetz number as follows:
\[
\Tr (\varphi \tei H^{\hullet}) := \sum_i (-1)^i \Tr (\varphi \tei H^i) \; .
\]
The following theorem is a reformulation of results in \cite{D1}. They are based on algebraic number theory and in particular on class field theory.

\begin{theorem}
  \label{t20}
Let $(\sigma , e)$ be an endomorphism of $(\eX , F)$ as above, $l \neq 2$. The cohomologies $H^{\hullet}_c (\eX , F)$ and $H^{\hullet} (\eX , F)$ are finite dimensional and we have:
\begin{equation}
  \label{eq:42}
  \Tr ((\sigma , e)^* \tei H^{\hullet}_c (\eX , F)) = 0 \; .
\end{equation}
Furthermore we have
\begin{equation}
  \label{eq:43}
  \Tr ((\sigma , e)^* \tei H^{\hullet} (\oeX , j_* F)) = \Tr ((\sigma , e)^* \tei H^{\hullet} (\eX^{\ann}_{\R} , \alpha^* F)) \; .
\end{equation}
If in addition $\eX$ is generically smooth and the fixed points $x$ of $\sigma$ on $\eX_{\infty}$ are non-degenerate in the sense that $\det (1 - T_x\sigma \tei T_x (\eX \otimes \R))$ is non-zero then we have:
\begin{equation}
  \label{eq:44}
  \Tr ((\sigma , e)^* \tei H^{\hullet} (\oeX , j_* F)) = \sum_{x \in \eX_{\infty} \atop \sigma x = x} \Tr (e_x \tei (\alpha^* F)_x) \varepsilon_x (\sigma) \; .
\end{equation}
Here
\[
\varepsilon_x (\sigma) = \sgn \det (1 - T_x \sigma \tei T_x (\eX \otimes \R)) \; .
\]
\end{theorem}

The following assertion is an immediate consequence of the theorem.
\begin{cor}
  \label{t21}
 Let $\sigma$ be an endomorphism of a scheme $\eX$ proper over $\spec \Z$. Then we have for any $l \neq 2$:
\begin{equation}
  \label{eq:45}
  \Tr (\sigma^* \tei H^{\hullet}_c (\eX , \Q_l)) = 0 
\end{equation}
and
\begin{equation}
  \label{eq:46}
  \chi (H^{\hullet} (\oeX , j_* \Q_l)) = \chi (H^{\hullet} (\eX^{\ann}_{\R} , \Q_l)) \; .
\end{equation}
If in addition $\eX$ has a smooth generic fibre and the fixed points of $\sigma$ on $\eX_{\infty}$ are non-degenerate we have:
\begin{equation}
  \label{eq:47}
  \Tr (\sigma^* \tei H^{\hullet} (\oeX , j_* \Q_l)) = \sum_{x \in \eX_{\infty} \atop \sigma x = x} \varepsilon_x (\sigma) \; .
\end{equation}
\end{cor}

\begin{punkt}
\rm We now explain analogies and an interesting difference with the case of varieties over finite fields. For a variety $X / \F_q$ with an endomorphism $\sigma$ of $X$ over $\F_q$ and a constructible $\Q_l$-sheaf $F$ on $X$ with an endomorphism $e : \sigma^* F \to F$ we have
\begin{equation}
  \label{eq:48}
  \Tr ((\sigma , e)^* \tei H^{\hullet}_c (X,F)) = 0
\end{equation}
and
\begin{equation}
  \label{eq:49}
  \Tr ((\sigma , e)^* \tei H^{\hullet} (X,F)) = 0 \; .
\end{equation}
Contrary to the number field case these assertions are easy to prove. E.g. for the first one, set $\oX = X \otimes \overline{\F}_q , \oF = F \, |_{\oX}$. The Hochschild--Serre spectral sequence degenerates into short exact sequences:
\[
0 \longrightarrow H^1 (\F_q , H^{n-1}_c (\oX , \oF)) \longrightarrow H^n_c (X , F) \longrightarrow H^0 (\F_q , H^n_c (\oX , \oF)) \longrightarrow 0 \; .
\]
Moreover for every $G_{\F_q}$-module $M$ there is an exact sequence:
\[
0 \longrightarrow H^0 (\F_q , M) \longrightarrow M \xrightarrow{1-\varphi} M \longrightarrow H^1 (\F_q , M) \longrightarrow 0 
\]
where $\varphi$ is a generator of $G_{\F_q}$. This implies (\ref{eq:48}) and (\ref{eq:49}) follows similarly since all groups are known to be finite dimensional.
\end{punkt}

Now in the number field case, according to Theorem \ref{t20} (\ref{eq:42}) the analogue of (\ref{eq:48}) is valid. In particular, for any finite set $S$ of prime ideals in $\eo_k$ we have
\begin{equation}
  \label{eq:50}
  \chi (H^{\hullet}_c (\spec \eo_{K,S} , \Q_l)) = 0 \; .
\end{equation}
The analogue of (\ref{eq:49}) however is not valid in general: From Theorem \ref{t20} (\ref{eq:43}) it follows for example that
\begin{equation}
  \label{eq:51}
  \chi (H^{\hullet} (\overline{\spec \eo_K} , j_* \Q_l) = \; \mbox{number of infinite places of} \; K \; .
\end{equation}
The fact, that in (\ref{eq:50}) the Euler characteristic is unchanged if we increase $S$ i.e. take out more finite places $\ep$ follows directly from the fact that
\[
\chi (H^{\hullet} (\spec (\eo_K / \ep) , \Q_l)) = 0 \; .
\]
The topological intuition behind this equation is that finite primes are like circles (whose Euler characteristic also vanishes). The difference between (\ref{eq:50}) and (\ref{eq:51}) comes from the different nature of the infinite places. Cohomologically the complex places behave like points and the real places behave like the ``quotient of a point by $G_{\R}$''. 

\begin{punkt}
  \rm In this section we prove a formula for Lefschetz numbers of finite order automorphisms of dynamical systems. We also allow certain constructible sheaves as coefficients.

Let us consider a complete flow $\phi^t$ on a compact manifold $X$. Assume that for every fixed point $x$ there is some $\delta_x > 0$ such that
\[
\det (1 - T_x \phi^t \tei T_x X) \neq 0 \quad \mbox{for} \quad 0 < t < \delta_x \; .
\]
This condition is weaker than non-degeneracy. Still it implies that the fixed points are isolated and hence finite in number.

We also assume that the lengths of the closed orbits are bounded below by some $\varepsilon > 0$.

Let $\sigma$ be an automorphism of finite order of $X$ which commutes with all $\phi^t$ i.e. an automorphism of $(X , \phi)$. Consider a constructible sheaf of $\Q$\,- or $\R$-vector spaces on $X$ together with an endomorphism $e : \sigma^{-1} F \to F$. We also assume that there is an action $\psi^t$ over $\phi^t$ i.e. isomorphisms $\psi^t : (\phi^t)^{-1} F \to F$ for all $t$ satisfying the relations
\[
\psi^0 = \id \quad \mbox{and} \quad \psi^{t_1 + t_2} = \psi^{t_2} \verk (\phi^{t_2})^{-1} (\psi^{t_1}) \quad \mbox{for all} \; t_1 , t_2 \in \R \; .
\]
For constant $F$ there is a canonical action $\psi$.\\
Note that $\sigma$ and $e$ determine an endomorphism $(\sigma , e)^*$ of $H^{\hullet} (X, F)$.\\
Let us say that a point $x \in X$ is $\phi$-fix if $\phi^t (x) = x$ for all $t \in \R$. 

Then we have the following formula:
\end{punkt}

\begin{theorem} \label{t24}
  $\dis \Tr ((\sigma , e)^* \tei H^{\hullet} (X, F)) = \sum_{x\in X , \phi- \mathrm{fix} \atop \sigma x = x} \Tr (e_x \tei F_x) \varepsilon_x (\sigma)$\\
where
\[
\varepsilon_x (\sigma) = \lim_{t \to 0 \atop t > 0} \sgn \det (1 - T_x (\phi^t \sigma ) \tei T_x X)\; .
\]
\end{theorem}

\begin{proof}
  Assume $N \ge 1$ is such that $\sigma^N = 1$ and fix some $s$ with \\
$0 < s < N^{-1} \min_{x\in X , \phi - \mathrm{fix}} (\varepsilon , \delta_x)$. If $x$ is any point of $X$ with $(\phi^s \sigma) (x) = x$, then $(\phi^{Ns} \sigma^N) (x) = x$ and hence $\phi^{Ns} (x) = x$. If $x$ lay on a periodic orbit $\gamma$ then $l (\gamma) \le Ns$ and hence $\varepsilon \le Ns$ contrary to the choice of $s$. Thus $x$ is a fixed point of $\phi$. Because of $\phi^s (x) = x$ we have $\sigma x = x$ as well. If $1$ is an eigenvalue of $T_x (\phi^t \sigma)$ then $1$ is an eigenvalue of $(T_x (\phi^t \sigma))^N = T_x \phi^{Nt}$ as well. Thus $Nt \ge \delta_x$. By assumption on $s$ we have $Ns < \delta_x$ and therefore
\[
\det (1 - T_x (\phi^s \sigma) \tei T_x X) \neq 0 \; .
\]
In conclusion: the fixed points of the automorphism $\phi^s \sigma$ coincide with those fixed points of the flow $\phi$ which are also kept fixed by $\sigma$. They are all non-degenerate.

Consider the morphisms for $t \in \R$
\[
e_t = e \verk \sigma^{-1} (\psi^t) : (\phi^t \sigma)^{-1} F \longrightarrow F \; .
\]
The Lefschetz fixed point formula for the endomorphism $(\phi^s \sigma , e_s)$ of $(X , F)$ for a fixed $s$ as above now gives the formula:
\[
\Tr ((\phi^s \sigma , e_s)^* \tei H^{\hullet} (X ,F)) = \sum_{x \in X ,  \phi-\mathrm{fix} \atop \sigma x = x } \Tr ((e_s)_x \tei F_x) \sgn \det (1 - T_x (\phi^s \sigma) \tei T_x X) \; .
\]
The left hand side is defined for all $s$ in $\R$ and by homotopy invariance of cohomology it is independent of $s$. Passing to the limit $s \to 0$ for positive $s$ in the formula thus gives the assertion.
\end{proof}

\begin{cor} \label{t25}
  Let $(X , \phi^t , \sigma)$ be as in theorem \ref{t24} and let $U \subset X$ be open, $\phi$- and $\sigma$-invariant. Assume that $X \ohne U$ is a compact submanifold of $X$. Consider a constructible sheaf $F$ of $\Q$- or $\R$-vector spaces on $U$ with an endomorphism $e : \sigma^{-1} F \to F$ and an action $\psi^t$ over $\phi^t \, |_U$. Then we have:
\[
\Tr ((\sigma , e)^* \tei H^{\hullet}_c (U , F)) = \sum_{x \in U , \phi-\mathrm{fix} \atop \sigma x = x} \Tr (e_x \tei F_x) \varepsilon_x (\sigma) \; .
\]
In particular
\[
\Tr ((\sigma , e)^* \tei H^{\hullet}_c (U ,F)) = 0
\]
if $\phi$ has no fixed points on $U$.
\end{cor}

\begin{proof}
  Apply \ref{t24} to $(X , j_! F)$, where $j : U \hookrightarrow X$ is the inclusion.
\end{proof}

\begin{punkt}
\rm  We now compare arithmetic and dynamic Lefschetz numbers. This section is of a heuristic nature: we assume that a functor $\Uh \mapsto (\Uh_{\dyn}, \phi^t)$ from flat algebraic schemes over $\spec \Z$ to dynamical systems exists, with properties as described before. We consider various cases.
\bigskip

\noindent {\sc 2.7.1} Let $\sigma$ be a finite order automorphism of a scheme $\eX$ proper and flat over $\spec \Z$. The induced automorphism $\sigma$ of $\eX_{\dyn}$ is an automorphism of $(\eX_{\dyn} , \phi^t)$ and hence commutes with each $\phi^t$. As we have seen, the phase space $\eX_{\dyn}$ cannot be a manifold. Still let us assume that the assertion of corollary \ref{t25} applies to
\[
(U , \phi^t , \sigma , F , e , \psi^t) = (\eX_{\dyn}, \phi^t , \sigma , \Q , \id , \id) \; .
\]
Then we find
\[
\Tr (\sigma^* \tei H^{\hullet}_c (\eX_{\dyn} , \Q)) = 0
\]
in accordance with formula (\ref{eq:45}) since $\phi^t$ should have no fixed points on $\eX_{\dyn}$.
\bigskip

\noindent {\sc 2.7.2} Now consider the case where in addition $\eX$ has a smooth generic fibre and the fixed points of $\sigma$ on $\eX_{\infty}$ are non-degenerate. Under assumptions as above, from corollary \ref{t25} applied to
\[
(U , \phi^t , \sigma, F, e , \psi^t) = (\oeX_{\dyn} , \phi^t , \sigma , j_* \Q , \id , \id)
\]
we would get a formula corresponding to (\ref{eq:47}):
\begin{equation}
  \label{eq:52}
  \Tr (\sigma^* \tei H^{\hullet} (\oeX_{\dyn} , j_* \Q)) = \sum_{x \in \eX_{\infty} \atop \sigma x = x} \varepsilon_x (\sigma) \; .
\end{equation}
Here we have used that the fixed point set of $\phi^t$ on $\oeX_{\dyn}$ should be $\eX (\C) / G_{\R} = \eX_{\infty}$. Namely, in the case of $\eX = \spec \eo_K$, this is the set of archimedean valuations of $K$. For general $\eX$ we only have the following argument: The set of closed points of $\eX$ over $p$ can be identified with the set $\eX (\OF_p) / \langle \Fr_p \rangle$ of Frobenius orbits on $\eX (\OF_p)$. Thus the set of closed orbits of $(\eX_{\dyn} , \phi^t)$ would be in bijection with the union of all $\eX (\OF_p) / \langle \Fr_p \rangle$. Correspondingly it looks natural to assume that the set of fixed points of $\phi$ would be in bijection with $\eX (\C) / \langle F_{\infty} \rangle = \eX_{\infty}$ where the infinite Frobenius $F_{\infty} : \eX (\C) \to \eX (\C)$ acts by complex conjugation. 

Actually the comparison of formulas (\ref{eq:46}) and (\ref{eq:52}) lends further credibility to this idea. 

The different definitions of $\varepsilon_x (\sigma)$ in (\ref{eq:46}) and (\ref{eq:52}) should agree. 
\end{punkt}

\noindent Mathematisches Institut\\
Einsteinstr. 62\\
48149 M\"unster, Germany\\
deninger@math.uni-muenster.de
\end{document}